\documentclass{article}
\usepackage[%
journal=AIHPD,
lang=american,
]{ems-journal}

\newcommand{\cE}{{\cal E}}

\newcommand{\cL}{{\cal L}}
\newcommand{\cM}{{\cal M}}
\newcommand{\cS}{{\cal S}}
\newcommand{\cX}{{\cal X}}
\newcommand{\cY}{{\cal Y}}
\newcommand{\ou}{{\overline{u}}}
\newcommand{\N}{{\mathds N}}
\newcommand{\R}{{\mathds R}}
\newcommand{\Z}{{\mathds Z}}
\newcommand{\Prob}{\text{Prob}}

\usepackage{dsfont}

\numberwithin{equation}{section}

\theoremstyle{plain}
\newtheorem{theorem}{Theorem}
\newtheorem{proposition}{Proposition}[section]

\begin{document}

\title{Random Domino Tilings and the Arctic Circle Theorem}
\titlemark{Random Domino Tilings and the Arctic Circle Theorem}

%

\emsauthor{1}{
	\givenname{William}
	\surname{Jockusch}
	\mrid{}
	\zblid{jockusch.william}
	\orcid{}}{W.~Jockusch}
\emsauthor{2}{
	\givenname{James}
	\surname{Propp}
	\mrid{}
	\zblid{propp.james-gary}
	\orcid{0000-0002-0913-3649}}{J.~Propp}
\emsauthor{3}{
	\givenname{Peter}
	\surname{Shor}
	\mrid{}
	\zblid{shor.peter-w}
	\orcid{}}{P.~Shor}

\Emsaffil{1}{
	\department{}
	\organisation{Wizards of the Coast}
	\rorid{}
	\address{1600 Lind Ave SW, Suite 400}
	\zip{98057}
	\city{Renton, WA}
	\country{U.S.A.}
	\affemail{jockusch@gmail.com}}
\Emsaffil{2}{
	\department{Department of Mathematics and Statistics}
	\organisation{University of Massachusetts Lowell}
	\rorid{03hamhx47}
	\address{1 University Avenue}
	\zip{01854}
	\city{Lowell, MA}
	\country{U.S.A.}
	\affemail{jamespropp@gmail.com}}
\Emsaffil{3}{
	\department{Mathematics}
	\organisation{Massachusetts Institute of Technology}
	\rorid{042nb2s44}
	\address{77 Massachusetts Avenue}
	\zip{02139}
	\city{Cambridge, MA}
	\country{U.S.A.}
	\affemail{shor@math.mit.edu}}

\classification{60C05}

\keywords{tiling, exclusion process}

In this article we study domino tilings of a family of finite regions
called Aztec diamonds.  Every such tiling determines a partition of the
Aztec diamond into five sub-regions; in the four outer sub-regions,
every tile lines up with nearby tiles, while in the fifth, central sub-region,
differently-oriented tiles co-exist side by side.  We show that when
$n$ is sufficiently large, the shape of the central sub-region
becomes arbitrarily close to a perfect circle of radius $n/\sqrt{2}$
for all but a negligible proportion of the tilings.
Our proof uses techniques from the theory of interacting particle systems.
In particular, we prove and make use of a classification of the stationary
behaviors of a totally asymmetric one-dimensional exclusion process
in discrete time.
\begin{abstract}
\ldots
\end{abstract}

\maketitle


\section{Introduction}
\label{sec:intro}

\begin{figure}
\begin{center}
\includegraphics[width=4.8in]{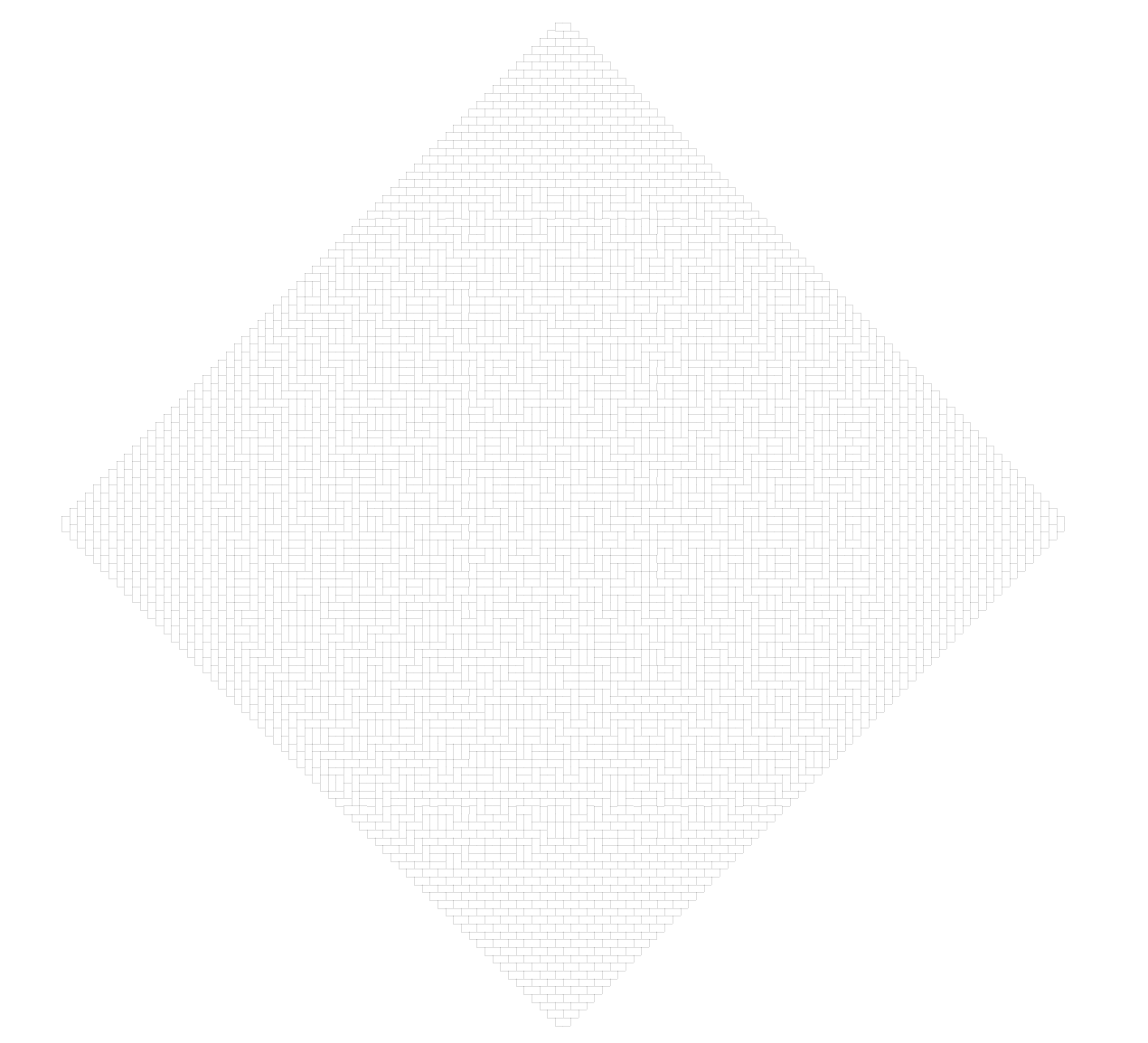}
\end{center}
\caption{A random domino tiling of the Aztec diamond of order 64.}
\label{fig1}
\end{figure}

This paper presents a complete and self-contained probabilistic proof 
of the original arctic circle theorem, written at the time of its discovery.

Figure~\ref{fig1} shows a domino tiling of the Aztec diamond of order 64.
In general, the {\bf \boldmath Aztec diamond of order $n$}
is a region composed of $2n(n+1)$ unit squares,
arranged as a stack of $2n$ centered rows of squares,
with the $k$th row having length min($2k,4n-2k+2$),
and a {\bf domino} is the union of any two of these unit squares
that share an edge.
It can be seen that the tiling shown in Figure~\ref{fig1} consists
of a roughly circular region in which tile-orientations are mixed,
surrounded by four regions in which the tiling exhibits
a repetitive ``brick-wall'' pattern.
In this article we will demonstrate that
the large-scale structure seen in the figure
is no coincidence but is in fact increasingly certain to occur
as one looks at random tilings of ever-larger Aztec diamonds.
This behavior is in sharp contrast with the behavior of
random domino tilings of $2n$-by-$2n$ squares,
which are statistically homogeneous 
unless one looks quite close to the boundary
\cite{burton}.

As a first step towards a precise statement of the main result,
it is important to note that the four brick-wall patterns
are genuinely different from one another.
One might not at first see the difference between
the brick-pattern seen at the top of the diamond
and the brick-pattern seen at the bottom,
but in fact they are out of phase with one another.
To see this, imagine coloring the squares alternately black and white.
Then, scanning from left to right,
rows in the upper half of the diamond begin with a square of one color
while rows in the lower half begin with a square of the other color.
Thus the two horizontal brick-wall patterns 
extend to distinct tilings of the entire plane;
in one, every domino has its left square colored black,
and in the other, every domino has its left square colored white.
A similar situation prevails for the vertical tiles.

The recognition of the existence of 
four rather than two distinct sorts of tiles
plays a key role in the formulation of the ``shuffling algorithm''
(described in section~\ref{sec:shuffling}).
Shuffling was introduced in \cite{elkies}
to prove that the Aztec diamond of order $n$
has exactly $2^{n(n+1)/2}$ domino tilings;
here we use shuffling to generate tilings uniformly at random
(Figure~\ref{fig1} was generated in precisely this way).

\begin{figure}
\begin{center}
\includegraphics[width=3.4in]{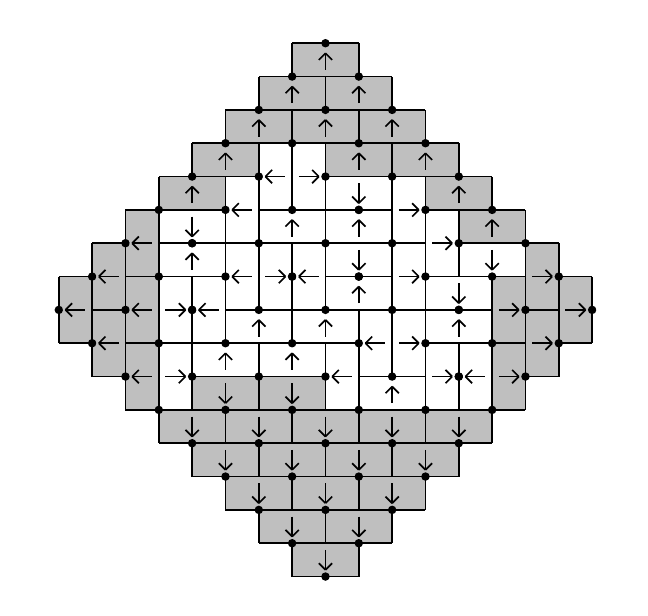}
\end{center}
\caption{Domino shuffling, before dominoes slide.}
\label{fig2}
\end{figure}

Figure~\ref{fig1} shows a domino tiling of the Aztec diamond of order 64.
To lay the groundwork for a description of shuffling,
consider the set of vertices 
in a paving of the Aztec diamond of order $n$ by unit squares.
Figure~\ref{fig2} shows a domino tiling
of the Aztec diamond of order 8.
As shown in the figure,
we mark with a dot the middle vertex 
along the upper (or ``northern'') border
of the Aztec diamond,
and we also mark with a dot
every vertex that can be reached from this vertex
by a lattice-path of even length.
Every domino then has a dot at the midpoint of 
exactly one of its four sides,
and we will assign the domino a {\bf heading}
(north, south, east, or west)
according to which of its sides 
(north, south, east, or west)
is dotted.
Thus a horizontal domino is either north-going or south-going
and a vertical domino is either east-going or west-going;
the reason for this terminology will become clear in section~\ref{sec:shuffling}.
Headings of dominoes are indicated by arrows in the figure.

Some of the dominoes in Figure~\ref{fig2} have been shaded,
namely, those that percolate to the boundary
by way of dominoes with the same heading.
For instance, a north-going domino has been shaded if and only if
it is possible to get from that domino to
a north-going domino that shares an edge with the boundary
of the Aztec diamond
by means of a (possibly trivial) sequence of north-going dominoes,
each of which shares an edge with the one before.

Observe that the only way a north-going domino
can share an edge with the boundary
is if it abuts one of the horizontal edges
in the upper half of the diamond,
and that the presence of such a domino leads to a cascade of 
north-going dominoes reaching all the way up to the top row.
It follows easily from this that
the shaded north-going dominoes
form a single connected block,
which we call the {\bf arctic region}.
(Similar remarks apply to the other three ``frozen'' regions.)
The unshaded region is called the {\bf temperate zone}.

Define the {\bf inscribed circle}
as the circle of radius $n/\sqrt{2}$
centered on the middle of the Aztec diamond.
We can now state the main result of this paper,
the Arctic Circle Theorem:

\begin{theorem} 
Fix $\epsilon > 0$.  Then for all sufficiently large $n$,
all but an $\epsilon$ fraction of the domino tilings
of the Aztec diamond of order $n$
will have a temperate zone whose boundary 
stays uniformly within distance $\epsilon n$
of the inscribed circle.
\end{theorem}

Note that this implies, for all $\epsilon > 0$,
that when $n$ is sufficiently large,
all but an $\epsilon$ fraction of the domino tilings
of the Aztec diamond of order $n$
will have a temperate zone whose
symmetric difference with the interior of the inscribed circle
has area less than $\epsilon n^2$.

To prove the Arctic Circle Theorem, we will use facts about
the totally asymmetric exclusion process on $\Z$
taking place in discrete time. 
This is a stochastic process whose states are
doubly-infinite sequences $(\dots,x_{-1},x_0,x_1,\dots)$ of 1's and 0's,
whose $i$th term signifies the presence or absence
of a particle at location $i$ in a 1-dimensional lattice;
at each discrete time-step, a particle that has a vacancy to its right
has a 50\% chance of moving one step to the right 
and a 50\% chance of staying put.
(Note that a particle with another particle immediately to its right
has no chance of moving until the location to its right becomes vacant.)
That is to say: If the system at time $n$ is in a particular state
$(\dots,x_{-1},x_0,x_1,\dots)$
and you want to advance the system to time $n+1$,
then find all $i$ such that $x_i=1$ and $x_{i+1}=0$
(these $i$'s are necessarily non-consecutive), and for each such $i$,
toss a fair coin (with independent coins for different values of $i$);
when it comes up heads, put $x'_i=0$ and $x'_{i+1}=1$,
and when it comes up tails leave $x'_i=1$ and $x'_{i+1}=0$.
Put $x'_i = x_i$ for all other $i$.
In this way one obtains a new doubly-infinite sequence
$(\dots,x'_{-1},x'_0,x'_1,\dots)$.
One repeats this process $x \stackrel{?}{\mapsto} x'$ 
for infinitely many iterations.
(Here the ``?'' is to remind us that $x'$
is a {\it stochastic} function of $x$.)

Now consider the initial condition $x^*$ with
\[
x_i^* = \left\{ \begin{array}{ll} 
	1 & \text{if $i \leq 0$,} \\
	0 & \text{if $i > 0$.}
	\end{array} \right.
\]
We are interested in how the exclusion process evolves
when it is in state $x^*$ at time 0.

It will be shown in section~\ref{sec:reduction}~that
the Arctic Circle Theorem can be reduced to the following assertion:

\begin{theorem}
Fix $\alpha \leq \beta$.
If one runs the exclusion process
starting from the state $x^*$,
then the number of particles in the interval $[\alpha n, \beta n]$
at time $n$, normalized by dividing by $n$,
converges in probability to $h(\alpha)-h(\beta)$, where
\[
h(u) = \left\{ \begin{array}{ll}
	-u & 
		\text{for $u<-\frac12$,} \\
	{\scriptstyle \frac{1-u}{2} - \frac12 \sqrt{\frac12-u^2}} & 
		\text{for $-\frac12 \leq u \leq \frac12$, and} \\
	0 & \text{for $u>\frac12$}. \end{array} \right.
\]
\end{theorem}

Theorem 2 is quite analogous to Rost's Theorem \cite{rost}
on the behavior of the totally asymmetric exclusion process
in {\it continuous} time.
The quantitative form of the result is different
(Rost obtains a parabolic arc where we obtain a circular arc),
but most of the methods are the same.

Note that the assertion for $\frac12 = \alpha \leq \beta$
follows immediately from the fact that
the position of the rightmost particle at time $n$
is binomially distributed with mean $n/2$ and variance $n/4$.
The case of $\alpha \leq \beta = -\frac12$ is similar
(one looks at the leftmost vacancy
instead of the rightmost particle).
Hence all of the interest of Theorem 2
lies in the case $-\frac12 \leq \alpha \leq \beta \leq \frac12$.

To prove Theorem 2, we will need to establish some general facts
about the exclusion process.
Let $\cX$ denote the set of sequences $x = (\dots,x_{-1},x_0,x_1,\dots)$
with terms in $\{0,1\}$,
and let 
$\cM(\cX)$ denote the set of probability measures on $\cX$,
relative to the usual Borel $\sigma$-algebra on $\cX$
generated by the cylinder sets $\{x \in \cX:
x_{i_1} = b_1,\dots, x_{i_k} = b_k\}$
($i_1,\dots,i_k \in \Z$, $b_1,\dots,b_k \in \{0,1\}$).
The law governing the exclusion process, run for one time-step,
gives a stochastic function from $\cX$ to $\cX$,
sending $x$ to some $x'$.
If $\mu \in \cM(\cX)$ has the property that
$\Prob_{\mu(x)} (x' \in A) = \Prob_{\mu(x)} (x \in A) = \mu(A)$
for all Borel sets $A \subset \cX$,
then we say $\mu$ is {\bf stationary} under the dynamics.
Equivalently, if we define the time-evolution map
$F: \cM(\cX) \rightarrow \cM(\cX)$
by the formula $\Prob_{F(\mu)} (x \in A) = \Prob_{\mu} (x' \in A)$, 
then $\mu$ is stationary
if and only if it is a fixed point of $F$.

The shift-operator $k \mapsto k+1$ on $\Z$
gives rise to a shift-map $T$ on $\cX$,
with $(Tx)_i = x_{i+1}$ for all $i$;
this in turn gives rise to a shift-operation $T^*$ on $\cM(\cX)$,
with $(T^* \mu)(A) = \mu(T^{-1} A) = \Prob_{\mu(x)}(Tx \in A)$.
If a measure $\mu$ satisfies $T^* \mu = \mu$,
we say that $\mu$ is {\bf translation-invariant}
(or {\bf shift-invariant}).

If for all Borel sets $S$
and all reals $\alpha$ with $0 \leq \alpha \leq 1$
we define 
\[
(\alpha \mu_1 + (1-\alpha) \mu_2)(S) = 
\alpha \mu_1 (S) + (1-\alpha) \mu_2 (S)
\]
and if we decree that $\mu_n \rightarrow \nu$
if and only if $\mu_n(A) \rightarrow \nu(A)$ for all cylinder sets $A$
(this is often called the {\bf weak topology}), 
and if we make use of the fact that a probability measure in $\cM(\cX)$
is determined by the measure it assigns to cylinder sets
(see for instance Theorem 3.1 in \cite{billingsley}),
then we can give $\cM(\cX)$ an affine structure and a topology.
Under these definitions,
the set of stationary, translation-invariant probability measures 
in $\cM(\cX)$ is a compact convex set.
Thus, every element of it can be written as a convex combination
of extremal elements.

Given $0 \leq d \leq 1$, put $s=\sqrt{d^2+(1-d)^2}$, and let
\begin{align*}
p(0) &= 1-d, \\
p(1) &= d, \\
q(0,0) &= \frac{s-d}{1-d}, \\
q(0,1) &= \frac{1-s}{1-d}, \\
q(1,0) &= \frac{1-s}{d} \text{, and}\\
q(1,1) &= \frac{s-(1-d)}{d},
\end{align*}
with values at the endpoints $d=0$ and $d=1$ being given by continuity.
Let $\mu_d$ be the unique translation-invariant probability measure on $\cX$
such that 
\begin{equation*}
\mu_d (\{x \in \cX:\ x_0 = b_0,\: x_1 = b_1,\: x_2 = b_2,\: \dots, 
x_{k-1} = b_{k-1},\: x_k = b_k \})
\end{equation*}
is equal to 
\begin{equation*}
p(b_0) q(b_0,b_1) q(b_1,b_2) \cdots q(b_{k-1},b_k)
\end{equation*}
for all $k \geq 1$ and all $b_0,\dots,b_k$ in $\{0,1\}$.
(Later we will frequently write $p(0)$ as $p_0$, $q(0,0)$ as $q_{00}$, etc.)
Note that this measure has the Markov property:
any event that is measurable with respect to $\dots,X_{-2},X_{-1}$
is conditionally independent of
any event that is measurable with respect to $X_1,X_2,\dots$ given $X_0$.

\begin{theorem}
The translation-invariant, stationary probability measures on $\cX$
are precisely those that can be expressed as
convex combinations of the measures $\mu_d$,
$0 \leq d \leq 1$.
\end{theorem}

Theorem 3 will be of independent interest
to those who study particle system models.
The proof is similar to the proof
that was devised by Rost for use in the continuous-time setting;
the only non-straightforward trick that
the discrete-time version of the proof requires 
is the coupling argument presented in section~\ref{sec:steady}.

The rest of this article is organized as follows.
In section~\ref{sec:shuffling},
we describe a combinatorial operation on tilings,
called shuffling,
that plays a fundamental role in the proof.
In section~\ref{sec:reduction},
we show that the Arctic Circle Theorem (Theorem 1)
reduces to our assertion (Theorem 2)
about the behavior of the discrete-time exclusion process
when the initial state is $x^*$.
In section~\ref{sec:steady}, we classify
the stationary, translation-invariant measures in $\cM(\cX)$
(Theorem 3).
In sections~\ref{sec:monotonicity}~and~\ref{sec:profile}~we show
that Theorem 3 implies Theorem 2,
which in turn causes Theorem 1 to topple into place.
In section~\ref{sec:bias},
we modify the notion of shuffling to allow for
the study of random tilings in which
one sort of domino (horizontal or vertical) is favored over the other,
and we explore the extent to which the methods used in the unbiased case
still allow one to determine the asymptotic shape of the temperate zone.
Finally, in section~\ref{sec:conclusion}, we offer some comments
that may help the reader to understand why the existence
of spontaneous large-scale structure in random tilings
is not as surprising as it might at first seem.

Since the appearance of the present work, 
the arctic circle phenomenon has been placed 
in a broader mathematical context and refined in several directions. 
Alternative proofs based on height functions and variational principles 
were developed by Cohn, Elkies, and Propp \cite{CEP},
providing a conceptual explanation of limit shapes for domino tilings,
and the work of Cohn, Kenyon and Propp \cite{CKP}
situated the result within a more flexible variational framework.
Determinantal and integrable-systems methods, 
drawing on Kasteleyn's determinantal approach \cite{kasteleyn},
have led to sharp asymptotic results 
for local statistics and fluctuations near the arctic boundary; 
in particular, Johansson \cite{johansson} showed that boundary fluctuations 
occur on the $n^{1/3}$ scale and converge to the Airy process. 
More generally, analogous limit-shape phenomena have been established 
for a wide class of tiling and dimer models, 
including boxed plane partitions and periodic planar dimers, 
culminating in general results 
of Kenyon, Okounkov, and Sheffield \cite{KOS}
describing frozen boundaries as algebraic curves. 
These later developments complement rather than replace the present argument, 
which gives a direct and elementary derivation 
of the arctic circle theorem for Aztec diamonds.

We emphasize that the goal of the present paper
is not to survey subsequent developments,
but to give a complete and self-contained
probabilistic proof of the original result
using the method of particle systems.

\section{Shuffling}
\label{sec:shuffling}

For the convenience of the reader, we restate (without proof)
the details of the shuffling algorithm introduced in \cite{elkies}.
We note in passing that a cleaner argument for the validity
of the algorithm was found later by the author, in reformulating
domino-shuffling in terms of perfect matchings; in particular, 
the case-analysis in the proof of the Lemma that appears in 
section 6 of \cite{elkies}. See \cite{propp} for details.

Domino shuffling is a stochastic procedure that turns
a domino tiling of the Aztec diamond of order $n-1$
into one of several domino tilings of the Aztec diamond of order $n$.
If one starts from the (empty) tiling of the Aztec diamond of order $0$
and applies shuffling $n$ times, the result is a uniformly random
domino tiling of the Aztec diamond of order $n$ --- 
that is, each of the $2^{n(n+1)/2}$ tilings has 
probability $2^{-n(n+1)/2}$ of being generated.

Consider an Aztec diamond of order $n-1$ tiled by dominoes,
where each domino has been assigned
a heading (north, south, east, or west)
as described in section~\ref{sec:intro}.
When two dominoes share a side of length 2,
they must be heading in opposite directions.
If the arrows point away from each other,
the two dominoes form a {\bf good block};
if the arrows point towards each other,
they form a {\bf bad block}.

\begin{figure}
\begin{center}
\includegraphics[width=3.4in]{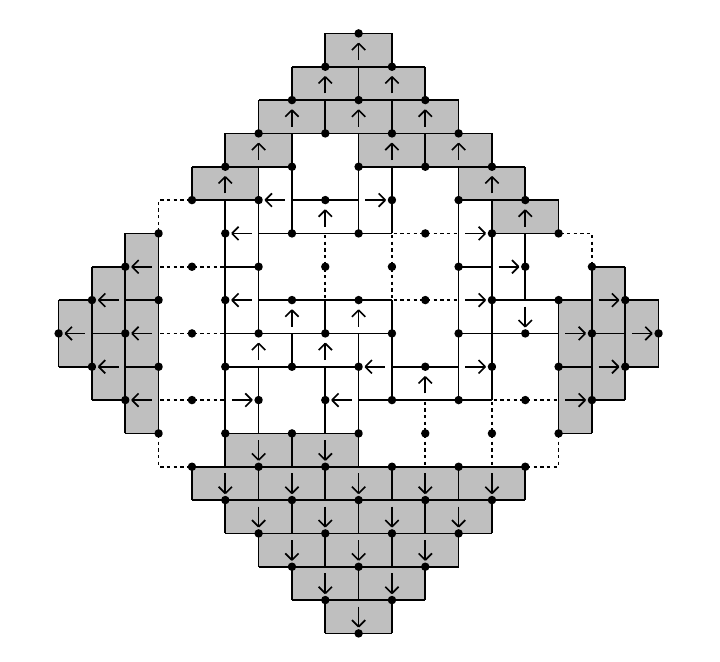}
\end{center}
\caption{Domino shuffling, after dominoes slide.}
\label{fig3}
\end{figure}

The shuffling procedure has three steps:
destruction, sliding, and creation.
Only the last of the three steps involves randomness.
In the destruction step, all the bad blocks
in the tiling of order $n-1$ are removed.
In the sliding step, each domino simultaneously moves
one step in the direction of its arrow.
After the sliding takes place,
some of the dominoes will no longer
be inside the Aztec diamond of order $n-1$,
but all will lie inside
the concentric Aztec diamond of order $n$.
It is shown in \cite{elkies} that
no dominoes overlap after destruction and sliding have taken place,
and that moreover the resulting configuration of dominoes,
viewed as a partial tiling of the Aztec diamond of order $n$
concentric with the original diamond of order $n-1$,
has a complement that can be tiled by 2-by-2 squares
in exactly one way.
In the creation step, we fill all these 2-by-2 holes
with good 2-by-2 blocks,
using a fair coin to decide whether the two tiles in any
particular hole should be horizontal or vertical.  The result is
a complete domino tiling of the Aztec diamond of order $n$.
Moreover, all the arrows in the new tiling
are properly assigned
in preparation for another round of shuffling.
(Note that the dotted and undotted vertices change places, as is fitting,
since we need the middle vertex 
on the upper border of the new Aztec diamond to be dotted
when we apply the shuffling algorithm.)
Figure~\ref{fig3} shows the results of applying
destruction and sliding
to the tiling shown in Figure~\ref{fig2};
the untiled portion of the Aztec diamond of order $n$
has been divided into 2-by-2 holes
in the only way possible.
When one performs creation,
the empty 2-by-2 holes get filled in
with good blocks,
and all the good blocks in the tiling of order $n$
arise in this way.

Shuffling can also be run in reverse:
given a tiling of the Aztec diamond of order $n$,
one can remove all the good 2-by-2 blocks,
slide each remaining tile one step in the direction
opposite to its arrow,
and then fill the 2-by-2 holes in the resulting partial tiling
of the Aztec diamond of order $n-1$
using bad blocks (each composed of either horizontal or vertical tiles).
We are not interested in the reverse algorithm per se;
rather, we will treat it as a way of seeing which tilings of the
smaller diamond can give rise to some specified tiling of the larger diamond
under the forward procedure.

To show that iterated shuffling yields 
the uniform distribution on the set of tilings,
it suffices to show that each round of shuffling
preserves uniformity.
For each tiling $T$ of order $n-1$,
let $b(T)$ be the number of bad 2-by-2 blocks
removed from the $T$ during the destruction step,
and for each tiling $T'$ of order $n$,
let $g(T')$ be the number of good 2-by-2 blocks
introduced when creating $T'$ during the creation step.
If $T$ and $T'$ are compatible in the sense that 
$T$ can give rise to $T'$ under one round of shuffling
(symbolically, $T \rightarrow T'$),
then by comparing the areas of the regions of order $n-1$ and $n$
we see we must have $g(T') - b(T) = n$.
Each tiling $T$ of the diamond of order $n-1$
can give rise to any of $2^{b(T)+n}$
different tilings $T'$ of the diamond of order $n$, 
each with probability $2^{-b(T)-n}$.  
Meanwhile, for each fixed $T'$, the number of preimages of $T'$ 
available to play the role of $T$ under one round of shuffling 
equals $2^{g(T')-n}$.
Hence, if we assume (for purposes of induction)
that each tiling $T$ of the Aztec diamond of order $n-1$ occurs with
probability $2^{-(n-1)n/2}$ after $n-1$ rounds of shuffling, then
after $n$ rounds, each tiling $T'$ of the Aztec diamond of order $n$
occurs with probability
\begin{align*}
\sum_{T: \ T \rightarrow T'} 2^{-(n-1)n/2} \ 2^{-b(T)-n} &= 
\sum_{T: \ T \rightarrow T'} 2^{-(n-1)n/2} \ 2^{-g(T')} \\
&= 2^{g(T')-n} \ 2^{-(n-1)n/2} \ 2^{-g(T')} \\
&= 2^{-n} \ 2^{-(n-1)n/2} \\
&= 2^{-n(n+1)/2}.
\end{align*}
The base case $n=0$ is trivial, so the claim follows.

\section{Reduction}
\label{sec:reduction}

To begin, let us check that the arctic region in
a domino tiling of an Aztec diamond
must have a fairly special sort of shape;
specifically, the centers of its constituent dominoes
must form a Ferrers diagram (rotated so that
its origin is at the north corner of the Aztec diamond).

\begin{figure}
\begin{center}
\includegraphics[width=3.4in]{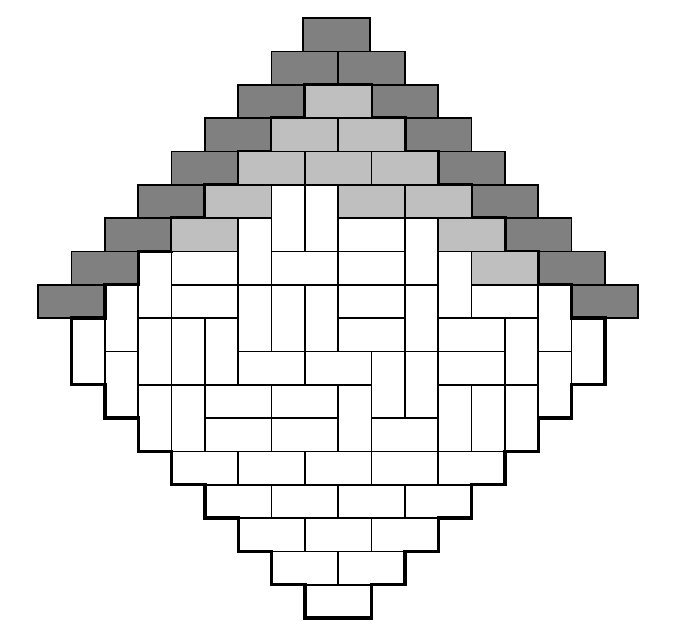}
\end{center}
\caption{The arctic region of a tiling,
along with external dominoes.}
\label{fig4}
\end{figure}

We will do this by giving an alternative characterization
of the arctic region.
Imagine (to make subsequent discussions simpler)
that we augment the picture by adding some extra dominoes
external to the Aztec diamond,
as shown in Figure~\ref{fig4}.
Specifically, we flank
each of the first $n-1$ rows of the Aztec diamond
by a horizontal domino on both the left and the right,
and we put two extra horizontal dominoes immediately above
the first row of the Aztec diamond,
with one more horizontal domino immediately above those two.
Call these $2n+1$ dominoes the {\bf external dominoes},
and say that they, together with the tiles inside the Aztec diamond,
constitute the {\bf augmented tiling}.
Now let us shade in the external dominoes,
and, in each successive row of the Aztec diamond
(starting with the first)
let us shade in every north-going domino
that is diagonally adjacent to
two already-shaded-in north-going dominoes in the preceding row,
and call those dominoes frozen. 
In Figure~\ref{fig4}, the external dominoes are given dark shading,
while the frozen dominoes that constitute the arctic region
are given light shading.

\begin{figure}
\begin{center}
\includegraphics[width=4.8in]{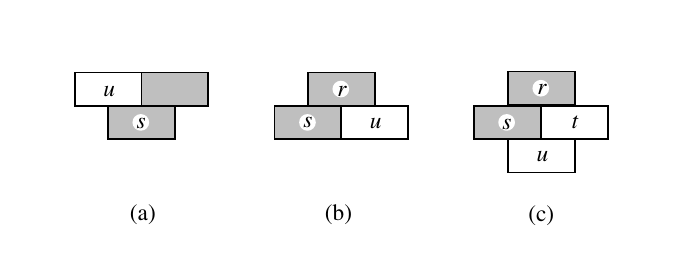}
\end{center}
\caption{Three cases.}
\label{fig5}
\end{figure}

We claim that no unshaded north-going domino
can share an edge with the shaded region.
For, suppose some (north-going) domino $s$ in the shaded region
shares an edge with some unshaded north-going domino inside the Aztec diamond.
Let $s$ be the northmost domino with this property,
and let $u$ be the unshaded north-going domino it borders.
For simplicity, we treat 
only the situation in which $s$ is one of the tiles
in our original tiling
(the argument is similar, but slightly simpler,
when $s$ is one of the external dominoes).
There are three cases, 
according to whether $u$ is in the row above $s$,
the same row as $s$, or the row below $s$.
The first case (shown in Figure~\ref{fig5}(a)) is easily disposed of;
since $s$ is shaded,
both dominoes covering it in the preceding row must be shaded,
contradicting the fact that $u$ is unshaded.
In the second case (shown in Figure~\ref{fig5}(b)),
the fact that $s$ is shaded
implies that the domino $r$ that covers both $s$ and $u$
must be shaded.
But then $s$ is not the northernmost shaded domino
in contact with an unshaded north-going domino
(as $r$ now has this property too).
In the third case (shown in Figure~\ref{fig5}(c)),
the fact that $s$ is shaded
implies that the domino $r$ must be shaded.
But then domino $t$ must be in the augmented tiling as well.
Since $u$ is unshaded, and $s$ is shaded,
$t$ must be unshaded;
but then $s$ is not the northernmost shaded domino
in contact with an unshaded north-going domino
(as $r$ now has this property too).

It follows that the shaded region
does not share an edge with any unshaded north-going dominoes.
(One can also see that 
it does not share a corner with any unshaded north-going dominoes either,
since two north-going dominoes never meet corner-to-corner.)
Therefore the part of the shaded region
that lies inside the Aztec diamond
coincides with the arctic region,
defined earlier as the union of the north-going dominoes
that percolate to the boundary.

The shaded region in the augmented tiling has the property that
each shaded domino in the original tiling 
is covered by two shaded dominoes
in the row above it,
each of which belongs to the augmented tiling.
Assign coordinates to the centers of the north-going dominoes
in accordance with a rotated rectangular coordinate system,
so that the northernmost external domino has center $(0,0)$
and the two below it have centers $(1,0)$ and $(0,1)$
(going from left to right).
Then the centers of the shaded dominoes are represented 
by a subset $S$ of $[0,n] \times [0,n]$
that contains $(0,k)$ and $(k,0)$ for all $0 \leq k \leq n$
and has the property that for all $(i,j) \in S$
with $1 \leq i,j \leq n$, both $(i-1,j)$ and $(i,j-1)$ are in $S$.
But such sets $S$ are precisely the Ferrers diagrams of partitions
with at most $n$ parts and with largest part at most $n$.

Now let us see how the arctic region changes under shuffling.
No domino in the arctic region belongs to a bad block.
Thus, north-going dominoes in the arctic region 
slide en masse to form a neighborhood of the ``north pole''
in ever-larger Aztec diamonds,
and remain part of the arctic region forever.
The question is, how do new dominoes get added to the arctic region over time?
No domino can join the arctic region by sliding;
the only way the arctic region can grow 
is by creation of good 2-by-2 blocks.

To apply shuffling to augmented tilings,
one simply applies shuffling to the part of the tiling
inside the Aztec diamond,
allows the external tiles to slide upward,
and creates two new external tiles
near the ``east and west poles''. 

\begin{figure}
\begin{center}
\includegraphics[width=4in]{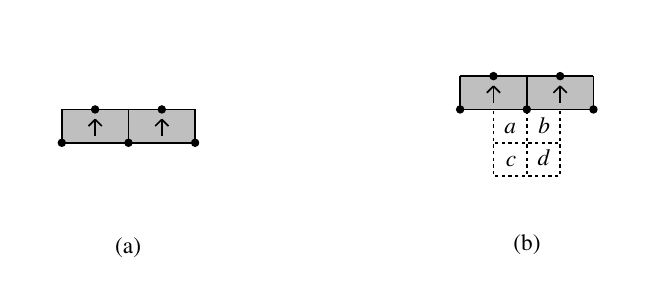}
\end{center}
\caption{How the arctic region grows.}
\label{fig6}
\end{figure}

Now consider two horizontally-adjacent north-going dominoes
in the augmented arctic region, as in Figure~\ref{fig6}(a).
Since we are interested in how the arctic region grows,
we assume we are at a growth-corner of the associated Ferrers diagram $S$;
that is, we assume that the two shaded dominoes
correspond to $(i-1,j) \in S$ and $(i,j-1) \in S$ with $(i,j) \not\in S$.
After sliding takes place,
the two dominoes
must be positioned as shown in Figure~\ref{fig6}(b).
Note that the cells marked $a$ and $b$
must be unoccupied at this point,
since east-, west-, and south-going dominoes
from the old (un-slid) tiling would not slide up that high.
Hence the cells marked $a$ and $b$ must
(in accordance with the results proved in \cite{elkies})
be part of an empty 2-by-2 block
(of the kind that becomes a good block when it is tiled),
namely, the 2-by-2 block shown in Figure~\ref{fig6}(b),
consisting of the cells marked $a$, $b$, $c$, and $d$.
Note that all such 2-by-2 holes are disjoint.
In the creation step, there is a $\frac12$ chance
that this empty block will be tiled by two horizontal dominoes,
in which case a new domino will get added to the arctic region,
and there is a $\frac12$ chance that the empty block
will be tiled by two vertical dominoes, in which case 
a new domino does not get added to the arctic region at that location.
Since the holes are disjoint, there is no interference
between the different ways different holes are filled;
each filling is determined by an independent coin-flip.
In terms of the Ferrers diagram, what is happening is that
whenever $(i,j) \not \in S$ with $(i-1,j) \in S$ and $(i,j-1) \in S$,
the node $(i,j)$ gets added to $S$ with probability $\frac12$.

Henceforth, let us disregard the external dominoes,
since they were only an aid to analysis
and play no further role in the proof.

It follows from the above argument,
via induction on $n$,
that after $n$ rounds of shuffling
the Ferrers diagram associated with the arctic region
contains only nodes $(i,j)$ with $i+j \leq n$;
that is, it sits inside
the Ferrers diagram of the partition $(n,n-1,\dots,1)$
(which is associated with the all-horizontals tiling).
Each Ferrers diagram that sits inside $(n,\dots,1)$ 
is in fact a possible shape of the arctic region,
though the different possibilities are not equally likely.

\begin{figure}
\begin{center}
\includegraphics[width=3in]{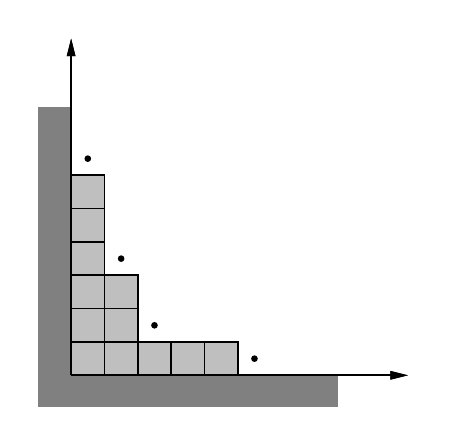}
\end{center}
\caption{The augmented Young diagram.}
\label{fig7}
\end{figure}

Let us represent the arctic region by a Young diagram
rather than by a Ferrers diagram.
That is, the $i$th domino-position in the $k$th row of the arctic region
now corresponds to the unit grid-square with lower-left corner $(k-i,i-1)$.
For convenience, imagine that the second, third, and fourth quadrants
are all adjoined to the Young diagram;
call this the {\bf augmented Young diagram}.
Figure~\ref{fig7} shows this augmented diagram, though it only includes 
the grid squares in the second, third, and fourth quadrants
that border the first quadrant.
Note that the Young diagram in the figure
corresponds to the arctic region shown in Figure~\ref{fig4}.
Call a square that is not in the augmented Young diagram
a ``growth-square''
if both its leftward and downward neighbors
are in the Young diagram.
Growth-squares in Figure~\ref{fig7} are marked by dots.

\begin{figure}
\begin{center}
\includegraphics[width=4in]{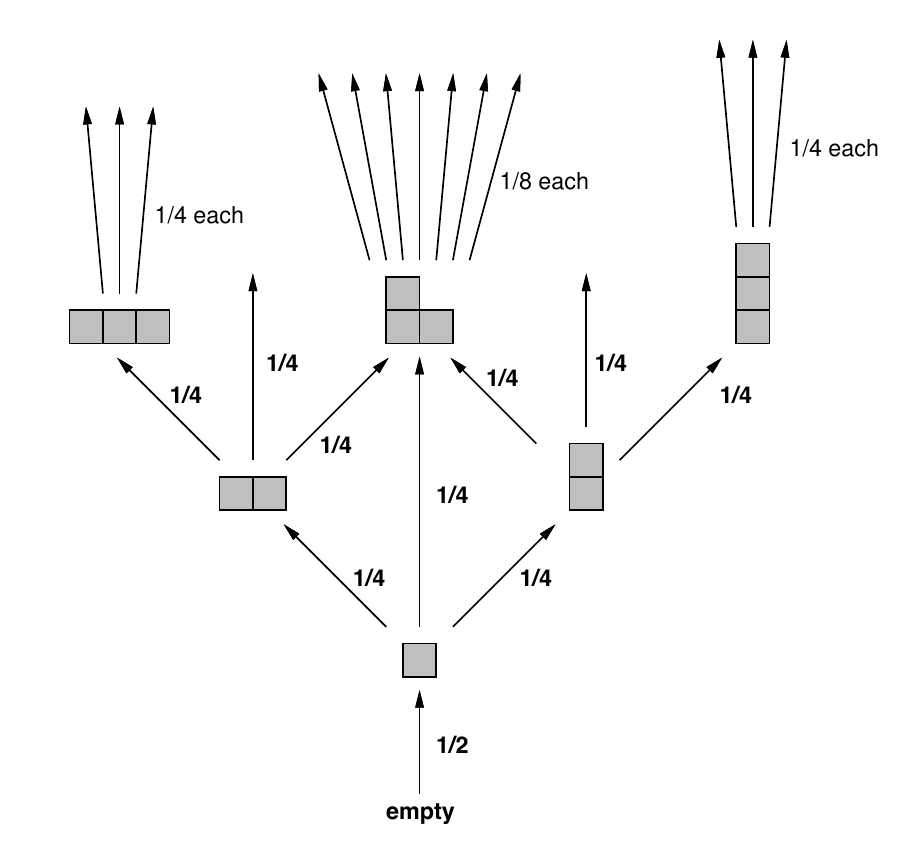}
\end{center}
\caption{Transition probabilities for growth.}
\label{fig8}
\end{figure}

The growth-process for Ferrers diagrams
(which represents the growth-process for possible shapes
of the arctic region) is tantamount to the growth process
for Young diagrams wherein, at each stage,
each of the growth-squares joins the growing diagram
independently with probability $\frac12$.
Figure~\ref{fig8} illustrates this
by giving the transition probabilities
for all Young diagrams with three or fewer boxes.
We have omitted the trivial transitions in which
a Young diagram gives rise to itself,
so the outgoing probabilities shown do not sum to 1.
Note also that the two two-box diagrams
admit transitions to four-box diagrams
that are not shown, though they are suggested by upward arrows.

A further transformation comes from looking at the boundary of
the augmented Young diagram.
This is a lattice-path with ``endpoints'' $(0,+\infty)$ and $(+\infty,0)$,
composed of unit steps downward and to the right.
The path initially runs straight
from $(0,+\infty)$ to $(0,0)$
and straight from $(0,0)$ to $(+\infty,0)$.
As the tiling undergoes iterated shuffling,
the lattice-path evolves 
in accordance with the rule that,
whenever a downward step is followed by a rightward step,
there is a probability of $\frac12$ that
``down-then-right'' will be replaced by ``right-then-down''
at the next stage,
with all such modifications being independent of one another.

We claim that after $n$ rounds of this evolution,
the lattice-path $L$ will,
with probability $1-o(1)$,
stay within distance $o(n)$ of a single particular path,
namely, the curvilinear path $P$ that goes 
from $(0,+\infty)$ to $(0,\frac n 2)$ along a straight line, 
from $(0,\frac n 2)$ to $(\frac n 2,0)$ 
along a quarter-circular arc
with its center at $(\frac n 2,\frac n 2)$,
and then from $(\frac n 2,0)$ to $(+\infty,0)$
along a straight line.
This is just a restatement of the main theorem.

To prove the restatement, it will suffice to prove that,
for all slopes $m$, $0<m<\infty$,
the intersection of $L$ with the line $y=mx$
will be within $o(n)$ of the intersection of $P$ with the line $y=mx$,
with probability $1-o(1)$.
For, by choosing a large but finite number of such slopes
that are suitably distributed 
between 0 and $\infty$,
we will be able to make sure that $P$ and $L$ stay uniformly close
merely by insuring that they are close at these check-points
(this argument makes use of the fact that $L$
can only go rightward or downward).

\begin{figure}
\begin{center}
\includegraphics[width=3.4in]{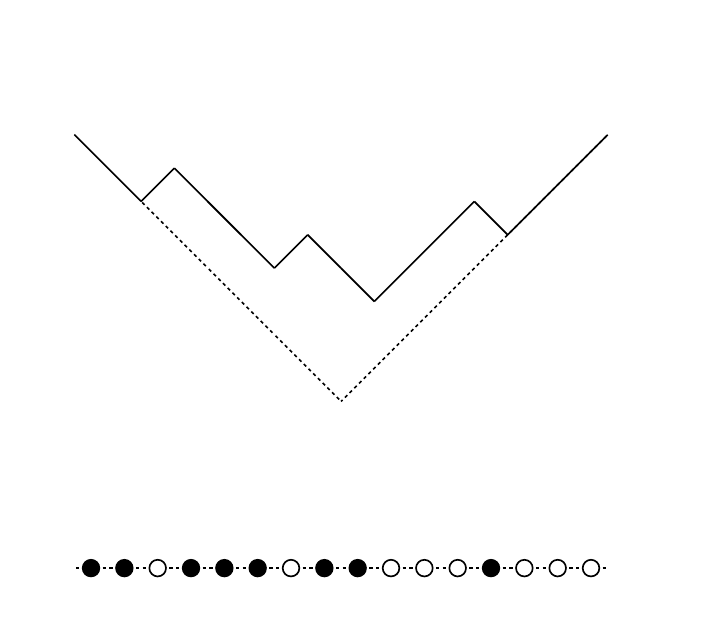}
\end{center}
\caption{Young diagram growth as a particle process.}
\label{fig9}
\end{figure}

As a final restatement, let us associate each lattice-path $L$ with 
a doubly infinite string $(\dots,x_{-1},x_0,x_1,\dots)$ of 0's and 1's.
Since $L$ consists of only rightward and downward steps,
each diagonal line $\{(x, y) : x-y = i\}$ 
is crossed exactly once by the monotone path $L$.
For $i \in {\Z}$ let $x_i$ be a 0 or a 1 
according to whether the unique step in $L$
from $L \cap \{(x,y) \mid x-y = i-1\}$
to $L \cap \{(x,y) \mid x-y = i\}$
goes rightward or downward.
Then the operation of changing a down-then-right jag
into a right-then-down jag in the lattice-path $L$ 
corresponds to the operation of 
replacing an occurrence of the substring ``10'' 
by an occurrence of the substring ``01'';
note that such substrings are necessarily disjoint,
and that each decision to replace 10 by 01 (or not),
just like each decision to modify a particular jag in $L$ (or not),
will be made independently of all other such decisions,
as determined by tosses of independent coins.
If we interpret a 1 in position $i$ as indicating the presence
of a particle at the $i$th site in a one-dimensional lattice,
and a 0 as indicating a vacancy there,
then this is just a rightward-jump event 
in an asymmetric exclusion process.
Figure~\ref{fig9} shows the lattice path of 
Figure~\ref{fig7} (tilted by 45 degrees for convenience)
and the associated configuration of the particle process.

Transition mechanisms of this kind were studied by Rost \cite{rost}.
However, Rost's dynamics occur in continuous time,
while ours take place in discrete time.
We will thus be able to avail ourselves of Rost's general approach,
but there will be some differences in the analysis.
In particular, the final results will not be the same:
we will obtain an arc of a circle
as the shape governing the system's asymptotic behavior,
where Rost obtains an arc of a parabola.

To conclude the link between Theorem 1 and 2,
we must show that the asymptotic circularity of the lattice-path $L$
described in the preceding paragraphs
is implied by the formula in Theorem 2.
To this end,
fix $0 < m < \infty$
and consider a lattice-point $(x,y)$ with $y = mx$.
The assertion that $(x,y)$ is on the lattice-path $L$ at time $n$
is equivalent to the assertion that in the corresponding state of
the particle process at time $n$,
there are $y$ particles to the right of position $x-y$.
If the density of particles is as given by Theorem 2,
then, putting $\beta=\frac12$ and $\alpha=(x-y)/n$,
we find that the number of particles to the right of position $x-y$
is $n(\frac{1-\alpha}{2} - \frac12\sqrt{\frac12-\alpha^2})$.
Equating this with $y$ and simplifying,
we obtain the relation
$(x-\frac{n}{2})^2+(y-\frac{n}{2})^2=\frac{n^2}{4}$,
describing a circle of radius $\frac{n}{2}$ 
centered on $(\frac{n}{2},\frac{n}{2})$.

\section{The Exclusion Process}
\label{sec:steady}

\subsection{Markov measures}
\label{subsec:markov}

The random variable $X_i$ is defined as the real-valued function on 
$\cX = \{0,1\}^\Z$
that sends $(\dots,x_{-1},x_0,x_1,\dots)$ to $x_i$.
We begin by proving that for all $0 \leq d \leq 1$
there is a unique stationary Markov measure $\mu=\mu_d$ such that
\begin{equation*}
p_1 = d
\end{equation*}
and
\begin{equation*}
q_{01} q_{10} = 2 q_{00} q_{11},
\end{equation*}
where
\begin{equation*}
p_i = \Prob_\mu [ X_0 = i ]
\end{equation*}
and
\begin{equation*}
q_{ij} = \Prob_\mu [ X_1 = j \mid X_0 = i ].
\end{equation*}
(The motivation behind this assertion
will emerge when we make use of 
the equation $q_{01} q_{10} = 2 q_{00} q_{11}$ below.)
For, a general shift-invariant Markov measure on $\cX$
is uniquely determined by $q_{00}$ and $q_{11}$
(since $q_{01} = 1-q_{00}$ and $q_{10} = 1-q_{11}$),
and the relation $q_{01} q_{10} = 2 q_{00} q_{11}$ yields
\begin{equation*}
(1-q_{00})(1-q_{11}) = 2q_{00}q_{11},
\end{equation*}
which can be solved for $q_{11}$:
\begin{equation*}
q_{11} = \frac{1-q_{00}}{1+q_{00}}.
\end{equation*}
So all that remains is to show that the condition $p_1 = d$
determines a unique value of $q_{11}$.
But note that
\begin{align*}
p_1 &= p_0 q_{01} + p_1 q_{11} \\
&= (1-p_1) q_{01} + p_1 q_{11} \\
&= q_{01} + p_1 (q_{11} - q_{01})
\end{align*}
so
\begin{align*}
p_1 &= \frac{q_{01}}{1-q_{11}+q_{01}} \\
&= \frac{1-q_{00}^2}{1+2q_{00}-q_{00}^2}.
\end{align*}
It's easy to check that this makes $p_1$
a strictly decreasing function of $q_{00}$
that goes from 1 to 0 as $q_{00}$ goes from 0 to 1.
Hence, for each $d$
with $0 \leq d \leq 1$
there is exactly one value of $q_{00}$
that yields $p_1 = d$, namely
\begin{equation*}
\frac{-d+\sqrt{d^2+(1-d)^2}}{1-d}
\end{equation*}
(when $d=1$ this expression is to be interpreted
as its limiting value, namely 0).
We will henceforth restrict our attention to the case $0<d<1$,
since $\mu_0$ and $\mu_1$ are trivial.
Note that $p_0 q_{01} = p_1 q_{10}$,
because the frequencies of the strings $0,1$ and $1,0$ must be equal.

Why this particular definition of $\mu_d$?
Imagine a finite version of the particle process introduced earlier,
in which the infinite line is replaced by a circle of $n$ sites,
$k$ of which are occupied by particles.
The update-rule is asymmetric as before:
particles may only advance ``to the right''
(counterclockwise, say).
The process can be modeled as a Markov chain
with $n \choose k$ states,
or as a random walk on a directed graph
with $n \choose k$ nodes.
If a configuration of the finite particle process
has $i$ instances of a particle with a vacancy to its right,
then the corresponding node of the directed graph
has $2^i$ outgoing arcs.
On the other hand, such a configuration also will have
$i$ instances of a vacancy with a particle to its right,
so the corresponding node of the directed graph
has $2^i$ incoming arcs.
It follows from this that a steady-state distribution
on the set of nodes
is given by a probability measure that assigns to each node
a probability proportional to $2^i$.
One can show that if one sends $k,n$ to infinity
with $k/n$ converging to some limit $d$,
then the statistics of these measures on circular words
converge weakly to the measure $\mu_d$.

Now we wish to show that $\mu = \mu_d$ is invariant
under the update-dynamics.
That is, we must show that
if we choose $x=(\dots,x_{-1},x_0,x_1,\dots) \in \cX$
in accordance with the distribution $\mu$
and then evolve the system through one time-step
to obtain a new doubly-infinite sequence $x' \in \cX$,
the probability of the event $x' \in B$
(denote this by $\mu'(B)$)
should be equal to 
the probability of the event $x \in B$
(which is just $\mu(B)$),
for all measurable events $B \subset \cX$.
To prove this,
it will suffice to prove $\mu'(B)=\mu(B)$ for cylinder sets $B$
of the form
\begin{equation*}
B = \{x \in \cX: x_0 = b_0 = 0,\ x_1 = b_1,\ x_2 = b_2, \dots,
\ x_{n-1} = b_{n-1}, \ x_n = b_n = 1\}
\end{equation*}
with $b_1,b_2,\dots,b_{n-1}$ in $\{0,1\}$ and with $n \geq 1$.
To see why this suffices, 
note first that the translation-invariance of $\mu_d$
and the fact that the update-dynamics commute with translation
guarantee that $\mu'$ will be translation-invariant.
Using this fact, combined with finite additivity, one can show that 
if $\mu'(B)=\mu(B)$ for all the special sets $B$ just described,
then $\mu'(B)=\mu(B)$ for all cylinder sets $B$
corresponding to bit strings that contain somewhere a 0 followed by a 1.
Indeed, let $C$ be the countable set consisting of all bit-strings
that nowhere contain a 0 followed by a 1; clearly $\mu(C)=0$.
Indeed, we must also have $\mu'(C) = 0$,
since the only way for $x'$ to belong to $C$
is if $x$ belongs to $C$.
Hence the cylinder sets containing at least one substring 01
generate the Borel $\sigma$-algebra up to 
a $\mu$-null set that is also $\mu'$-null.
Therefore, using countable additivity,
one can prove $\mu'(B)=\mu(B)$ for all cylinder sets
corresponding to finite bit-strings.
It follows from this that $\mu'$ agrees with $\mu$
on all cylinder sets
and hence on the entire measure-algebra of Borel sets.

Let $N(B) \geq 1$ be the number of occurrences of the substring 0,1
in the bit-string $b = (b_0,\dots,b_n)$.
There are $2^{N(B)}$ different bit-strings $a=(a_0,\dots,a_n)$
that can be formed by replacing none, some, or all
of these substrings by 1,0.
A point $x \in \cX$ can evolve in one time-step into a point $x' \in B$
if and only if $(x_1,\dots,x_n)$ is one of these $2^{N(B)}$ bit-strings.
Let $A = A(a)$ denote the set of $x \in \cX$
satisfying $x_0=a_0,\dots,x_n=a_n$,
and let $M(A)$ denote the number of occurrences
of $(1,0)$ in $(a_0,\dots,a_n)$,
Then the probability that an $x$ in $A$
will evolve in one time-step into an $x'$ in $B$
is equal to $2^{-M(A)}$
times two correction factors $r(A)$ and $s(A)$
associated with the beginning and end of the string $a$, respectively.

Specifically, if $a_n = 1$, then there are two ways 
in which the $n$th bit of $x$ could remain a 1,
depending on whether $x_{n+1}$ is 1 or 0.
In the former case (an event with probability $q_{11}$),
the $n$th bit must stay equal to 1, while in the latter case,
the $n$th bit will stay 1 with probability $\frac12$.
Thus the probability that the $n$th bit remains a 1
is $q_{11} + \frac12 q_{10}$. Since (as previously noted)
$q_{11} = (1-q_{00})/(1+q_{00})$ and $q_{10} = 1-q_{11}$,
we have $q_{10} = 1 - (1-q_{00})/(1+q_{00}) = 2 q_{00} / (1+q_{00})$
so that $q_{11} + \frac12 q_{10} = 
(1-q_{00})/(1+q_{00}) + q_{00} / (1+q_{00}) = 1/(1+q_{00})$.
Thus we define $s(A)$ to be $1/(1+q_{00})$ if $a_n = 1$
and to be 1 if $a_n = 0$.

Similarly, if $a_0 = 0$
the probability that the 0 remains a 0
can be shown to be $\frac12 (1+q_{00})$,
so we define $r(A)$ to be $\frac12 (1+q_{00})$ if $a_0 = 0$
and to be 1 if $a_0 = 1$. Then
\begin{equation*}
\mu'(B) = \sum_A \mu(A) \: 2^{-M(A)} r(A) s(A),
\end{equation*}
where the sum is over the $2^{N(B)}$ distinct cylinder sets $A$
associated with the bit-strings $a$ described above.

To prove that $\mu'(B) = \mu(B)$, it will suffice to prove that
each of the $2^{N(B)}$ summands is equal to $2^{-N(B)} \mu(B)$.
That is, we will show that every $A$
that contributes to the sum satisfies
$\mu(A)/\mu(B) = 2^{M(A)-N(B)} / r(A) s(A)$.
We do this by induction on the number of swaps
required to turn $B$ into $A$.
In the case of 0 swaps ($A=B$), the formula is true,
since $M(A)=N(B)-1$ and 
$r(A)s(A)=\frac12 (1+q_{00}) \cdot 1/(1+q_{00}) = \frac12$.
To get the induction step, we need to show that
\begin{equation}
\label{proportion}
\mu(\hat{A})/\mu(A) = 2^{M(\hat{A})-M(A)} 
r(A) s(A) / r(\hat{A}) s(\hat{A})
\end{equation}
when $A$ and $\hat{A}$ differ only in a single swap
(that is, $\hat{A}$ has 1,0
in two adjacent positions in which $A$ has 0,1).

The case $n=1$ is special.
In this circumstance, we have $a=(0,1)$ and $\hat{a}=(1,0)$,
with $\mu(\hat{A})/\mu(A) = 1$
on general principles.
On the right hand side of \eqref{proportion}
we get $M(\hat{A})-M(A)=1$,
$r(A) = \frac12 (1+q_{00})$,
$s(A) = 1 / (1+q_{00})$,
$r(\hat{A}) = 1$, and $s(\hat{A}) = 1$,
so the right hand side of (\ref{proportion}) is equal to 
\begin{equation*}
(2^1) ({\scriptstyle \frac12} (1+q_{00})) (1 / (1+q_{00})) = 1.
\end{equation*}
This verifies (\ref{proportion}) in the case $n=1$;
henceforth, we assume $n>1$.

Let us assume that the swap occurs in positions $i$ and $i+1$,
with $0 \leq i < i+1 \leq n$.

Suppose first that $i > 0$ and $i+1 < n$,
so that $r(\hat{A})=r(A)$
and $s(\hat{A})=s(A)$.
Then there are four cases to consider, 
according to the values of $a_{i-1}$ and $a_{i+2}$.
If $a_{i-1} = 0$ and $a_{i+2} = 0$,
then $M(\hat{A}) = M(A)$.
To find the ratio of $\mu(\hat{A})$ and $\mu(A)$,
write
\begin{equation*}
\mu(A) = p_{a_0} q_{a_0 a_1} \cdots q_{01} q_{10} q_{00} \dots q_{a_{n-1} a_n}
\end{equation*}
and
\begin{equation*}
\mu(\hat{A}) = 
p_{a_0} q_{a_0 a_1} \cdots q_{00} q_{01} q_{10} \dots q_{a_{n-1} a_n}.
\end{equation*}
The two products are just re-arrangements of one another,
so $\mu(\hat{A})/\mu(A)=1$.
If $a_{i-1} = 0$ and $a_{i+2} = 1$,
then $M(\hat{A}) = M(A)+1$
and $\mu(\hat{A})/\mu(A) = q_{01} q_{10} / q_{00} q_{11} = 2$.
If $a_{i-1} = 1$ and $a_{i+2} = 0$,
then $M(\hat{A}) = M(A)-1$
and $\mu(\hat{A})/\mu(A) = q_{11} q_{00} / q_{01} q_{10} = \frac12$.
Lastly,
if $a_{i-1} = 1$ and $a_{i+2} = 1$,
then $M(\hat{A}) = M(A)$
and $\mu(\hat{A})/\mu(A) = 1$.
In all four cases,
(\ref{proportion}) is verified.

When $i=0$,
there are two cases to consider,
according to the value of $a_2$.
If $a_2=0$, 
$\mu(\hat{A})/\mu(A) = p_1 q_{10} q_{00} / p_0 q_{01} q_{10}$.
Using the fact that $p_1 q_{10} = p_0 q_{01}$, this becomes
$q_{00} / q_{10}$, which equals $r(A)$.
Since $M(\hat{A}) = M(A)$
and $r(\hat{A}) = 1$ and $s(\hat{A}) = s(A)$,
(\ref{proportion}) holds.
If $a_2=1$, 
$\mu(\hat{A})/\mu(A) = 
p_1 q_{10} q_{01} / p_0 q_{01} q_{11} =
q_{01} / q_{11} = 1+q_{00} = 2r(A)$.
Since $M(\hat{A}) = 1+M(A)$
and $r(\hat{A}) = 1$ and $s(\hat{A}) = s(A)$,
(\ref{proportion}) holds.

When $i=n-1$,
there are two cases to consider,
according to the value of $a_{n-2}$.
If $a_{n-2}=0$, 
$\mu(\hat{A})/\mu(A) = q_{10} / q_{00} = 2/(1+q_{00}) = 2s(A)$.
Since $M(\hat{A}) = 1+M(A)$
and $s(\hat{A}) = 1$ and $r(\hat{A}) = r(A)$,
(\ref{proportion}) holds.
If $a_{n-2}=1$, 
$\mu(\hat{A})/\mu(A) = q_{11} / q_{01} = 1/(1+q_{00}) = s(A)$.
Since $M(\hat{A}) = M(A)$
and $s(\hat{A}) = 1$ and $r(\hat{A}) = r(A)$,
(\ref{proportion}) holds.

This completes the proof of
(\ref{proportion}),
which completes the proof that $\mu'(B) = \mu(B)$
for all our special cylinder sets $B$.
As remarked earlier,
this suffices to establish
that $\mu' = \mu$;
that is, $\mu = \mu_d$
is invariant under our evolution rules.

Note for later purposes that
$\mu(x_0=1,\ x_1=0)
=\mu(x_0=0,\ x_1=1)
=p_0 q_{01} = 1 - \sqrt{d^2+(1-d)^2}$.

\subsection{Uniqueness}
\label{subsec:uniqueness}

In the previous subsection
we found a one-parameter family of translation-invariant,
dynamically-stationary probability measures $\mu_d$.
As remarked earlier,
the translation-invariant, dynamically-stationary measures
form a compact convex subset 
of the compact space $\cM(\cX)$.
We will show that the only extremal points of this subset
are the measures $\mu_d$.
This will imply that every translation-invariant stationary measure
is a convex combination of the $\mu_d$'s.

Let $\mu$ be an extremal measure in 
the set of translation-invariant, dynamically-stationary measures,
with $\Prob_\mu [X_0 = 1] = d$.
We must show that $\mu = \mu_d$.
In preparation for a coupling argument,
let $\cX^{(1)} = \cX^{(2)} = \cX$
and let $\cM_{\mu,\mu_d}(\cX^{(1)} \times \cX^{(2)})$
be the space of translation-invariant probability measures on
$\cX^{(1)} \times \cX^{(2)}$
that project to $\mu$ on $\cX^{(1)}$
and $\mu_d$ on $\cX^{(2)}$.
Note that $\cM_{\mu,\mu_d}(\cX^{(1)} \times \cX^{(2)})$
is non-empty, since in particular $\mu \times \mu_d$ is in it.
Let
\begin{equation*}
\delta = \inf \big\{ \Prob_\pi [ X_0^{(1)} \neq X_0^{(2)} ] 
: \pi \in \cM_{\mu,\mu_d} (\cX^{(1)} \times \cX^{(2)}) \big\} .
\end{equation*}
$\cM_{\mu,\mu_d}$ is compact in its weak topology,
so the continuous function
$\pi \mapsto \Prob_\pi [ X_0^{(1)} \neq X_0^{(2)} ]$
achieves its minimum value at some particular $\pi$.
Henceforth, $\pi$ will denote just such a discrepancy-minimizing measure.
More specifically, since the discrepancy functional is affine,
the set of discrepancy-minimizing measures must be compact and convex, 
and invoking Krein-Milman we may choose $\pi$
to be an extremal point of this set.
Our strategy will be to show that if $\delta > 0$,
then we can find another joining $\pi'$ with smaller discrepancy,
contradicting the definition of $\delta$.

To construct $\pi'$ from $\pi$,
define dynamics on $\cX^{(1)} \times \cX^{(2)}$ as follows.
Fix $(x^{(1)}, x^{(2)}) \in \cX^{(1)} \times \cX^{(2)}$.
Say that the spatial indices $i$ and $i+1$ are {\bf linked}
(relative to $(x^{(1)}, x^{(2)})$)
if either of the doubly-infinite strings $x^{(1)}$, $x^{(2)}$
has a 1 in the $i$th position followed by a 0 in the $i+1$st position.
If positions $i$, $i+1$, $i+2$, \dots, $j-1$, and $j$ are pairwise linked
each to the next,
then we say that $i,i+1,\dots,j$ are in the same block.
In this way, $\Z$ is divided into blocks,
and $(x^{(1)}, x^{(2)})$ is divided into sub-words
(which we will also call blocks),
each of which consists of $k$ consecutive symbols from $x^{(1)}$
and the corresponding $k$ symbols from $x^{(2)}$, for some $k$.
(We might worry that there could be an infinite block,
but this would require an infinite stretch of values of $i$ for which
$(x^{(1)}_i,x^{(2)}_i)$ alternates between $(1,0)$ and $(0,1)$;
this in turn would imply that $x^{(2)}$
contains an infinitely long subword
in which 0's and 1's alternate,
and since $x^{(2)}$ is governed by $\mu_d$,
the chance of this happening is zero.)

An example will be helpful here:
if $(x^{(1)},x^{(2)})$ is
\begin{equation*}
\begin{array}{cccccccccccccc}
\dots & 0 & 0 & 1 & 1 & 0 & 1 & 1 & 0 & 1 & 1 & 0 & 1 & \dots \\
\dots & 0 & 0 & 1 & 0 & 1 & 0 & 0 & 0 & 0 & 1 & 0 & 0 & \dots
\end{array}
\end{equation*}
(with top and bottom rows corresponding to $x^{(1)}$ and $x^{(2)}$,
respectively),
then the blocks are
\[
\begin{array}{cccccccccccccccccccc}
\dots & 0 & \ & 0 & \ & 1 & 1 & 0 & 1 & \ & 1 & 0 & \ & 1 & \ & 1 & 0 & \ & 1 & \dots \\
\dots & 0 & \ & 0 & \ & 1 & 0 & 1 & 0 & \ & 0 & 0 & \ & 0 & \ & 1 & 0 & \ & 0 & \dots 
\end{array} .
\]
To describe the update-dynamics, 
one need only specify how each block is to be updated.

If a block is of length 1, there is nothing to do.  
If a block is of length 2, 
then one has no choice about how to couple the two update-processes, 
unless the block is of the form
\[
\begin{array}{ll}
	1 & 0 \\
	1 & 0.
\end{array}
\]
In this case, 
we decree that one should use
the {\it same} random bit 
to decide what to do with the upper substring (from $x^{(1)}$) 
and the lower substring (from $x^{(2)}$).  
In particular, with probability $\frac12$ the block stays 
the same, and with probability $\frac12$ it becomes
\[
\begin{array}{ll}
	0 & 1 \\
	0 & 1.
\end{array}
\]

For blocks of length 3 or more, which are of the form
\[
\begin{array}{ccccccc}
1 & 0 & 1 & 0 & \dots & 1 & 0 \\
? & 1 & 0 & 1 & \dots & 0 & ?
\end{array}
\]
(or minor variations thereof), 
one must do something slightly different: 
we decree that 
in deciding about the $i$th occurrence of ``1 0'' in the upper substring
and the $i$th occurrence of ``1 0'' in the lower substring,
one should use not the same bit but {\it complementary} bits.
That is to say,
one should convert
the $i$th occurrence of ``1 0'' in the $x^{(1)}$-row of the block to a ``0 1'' 
if and only if one leaves the $i$th occurrence of ``1 0'' 
in the $x^{(2)}$-row of the block alone.
In all other respects 
one's random choices are to be independent of one another.

Define a mismatch as 
a position in which the $x^{(1)}$ and $x^{(2)}$ words disagree.
If a block has length $n \geq 2$,
the number of mismatches goes from at least $n-2$
to at most $\lfloor (n-1)/2 \rfloor$.
This implies that in each finite block, 
the number of mismatches cannot increase.  
Indeed, in certain kinds of blocks 
(call them ``unstable blocks'')
the number of mismatches will go down with positive probability.  
The only blocks that do not have this property ---
the ``stable blocks,'' as we'll call them --- 
are the blocks of length 1 and the blocks
\[
\begin{array}{lllllllllllllll}
1 & 0 & \ & \ & 1 & 1 & \ & \ & 0 & 0 & \ & \ & 1 & 1 & 0 \\
1 & 0 & , & \ & 1 & 0 & , & \ & 1 & 0 & , & \ & 1 & 0 & 0 
\end{array}
\]
(and the blocks obtained from them 
by switching the roles of $x^{(1)}$ and $x^{(2)}$).

If we let $\pi_n$ be the probability measure
describing the outcome of applying the joint dynamics to $\pi$
for $n$ time-steps,
and we let $\Delta = \{(x^{(1)},x^{(2)}) \ : \ x_0^{(1)} \neq x_0^{(2)} \}$,
then (making use of the shift-invariance of the $\pi_n$'s)
we have $\pi(\Delta) \geq \pi_1 (\Delta) \geq \pi_2 (\Delta) \geq \dots\,$.
Since $\pi$ was chosen to minimize $\pi(\Delta) = \delta$,
we must in fact have $\pi(\Delta) = \pi_1(\Delta) = \pi_2(\Delta) = \dots\,$.
This means that unstable blocks have probability zero under $\pi$.

Now, we have
\begin{equation*}
\delta = \pi(\Delta) = 
\Prob_\pi [ X_0^{(1)} = 1, \ X_0^{(2)} = 0 ] +
\Prob_\pi [ X_0^{(1)} = 0, \ X_0^{(2)} = 1 ].
\end{equation*}
But
\begin{align*}
\Prob_\pi [ X_0^{(1)} = 1, \ X_0^{(2)} = 0 ]
&= \Prob_\pi [ X_0^{(1)} = 1 ] - 
   \Prob_\pi [ X_0^{(1)} = 1 , \ X_0^{(2)} = 1 ] \\
&= d - \Prob_\pi [ X_0^{(1)} = 1 , \ X_0^{(2)} = 1 ] \\
&= \Prob_\pi [ X_0^{(2)} = 1 ] - 
   \Prob_\pi [ X_0^{(1)} = 1 , \ X_0^{(2)} = 1 ] \\
&= \Prob_\pi [ X_0^{(1)} = 0 , \ X_0^{(2)} = 1 ] .
\end{align*}
Hence
$\Prob_\pi [ X_0^{(1)} = 1 , \ X_0^{(2)} = 0 ]
= \Prob_\pi [ X_0^{(1)} = 0 , \ X_0^{(2)} = 1 ]
= \delta/2$.
We now claim that this implies that some event of the form
$X_i^{(1)} = 1, \ X_i^{(2)} = 0, \ X_j^{(1)} = 0, \ X_j^{(2)} = 1$
must have positive probability.
For, were this not the case,
$\pi$ would assign probability 1 to
$\cY_1 \cup \cY_2$,
where $\cY_1 = \{ (x^{(1)},x^{(2)}): 
\ \text{$x_i^{(1)} \geq x_i^{(2)}$ for all $i$}\}$
and $\cY_2 = \{ (x^{(1)},x^{(2)}): 
\ \text{$x_i^{(1)} \leq x_i^{(2)}$ for all $i$}\}$.
Notice, however, that each of the two sets
is translation-invariant and is mapped into itself
by our evolution-rules;
if $\cY_1$ and $\cY_2$ each had positive probability,
then by restricting $\pi$ to each of them in turn,
we would get two distinct measures whose non-trivial weighted average was $\pi$,
contradicting the extremality of $\pi$.
Hence one of the two sets would have to have measure zero. 
But this contradicts the earlier-proved fact that
$\Prob_\pi [ X_0^{(1)} = 1 , \ X_0^{(2)} = 0 ]$ and
$\Prob_\pi [ X_0^{(1)} = 0 , \ X_0^{(2)} = 1 ]$
are both positive.
Hence, as asserted, there exist $i,j$ such that
$\Prob [ X_i^{(1)} = 1, \ X_i^{(2)} = 0, \ X_j^{(1)} = 0, \ X_j^{(2)} = 1]$
is positive.
Without loss of generality, assume $i<j$.

Consider now a two-rowed pattern of length $j-i+1$ of the form
\[
\begin{array}{lll}
	1 & \dots & 0 \\
	0 & \dots & 1 
\end{array}
\]
for some fixed (but unconstrained) values of the intervening bits,
occurring with positive probability under $\pi$.
We may assume without loss of generality that 
there are no discrepancies between the respective intervening bits,
for if this is not the case we can take $j-i$ smaller.
Let $A$ be the set of $(x^{(1)},x^{(2)}) \in \cX^{(1)} \times \cX^{(2)}$
with this pattern at positions $i$ through $j$.
By hypothesis, $\pi(A) > 0$.

Pick $(x^{(1)},x^{(2)}) \in A$ randomly,
in accordance with the measure $\pi$.
There is a positive probability that in some finite number of steps
the intervening bits will sort themselves so that 
the 1's are as far to the right as possible while
the 0's are as far to the left as possible.
At this point, assuming $j-i \geq 2$,
either the 1 on the left boundary can move to the right
(while the other positions are unaffected)
or the 1 on the right boundary can move to the left
(while the other positions are unaffected).
This reduces the distance between the discrepancies by 1.
Iterating this until we get down to $j-i=1$,
we see that there exists some $n$ and some $i$
for which
$\Prob_{\pi_n}
[ X_{i}^{(1)} = 1, \ X_{i}^{(2)} = 0, \ X_{i+1}^{(1)} = 0, \ X_{i+1}^{(2)} = 1 ]
> 0$.
However, such an $i$ would belong either to
a block of the form
\[
\begin{array}{cc}
1 & 0 \\
0 & 1
\end{array}
\]
or to a block of length $\geq 3$,
and this unstable block would lead to
$\Prob_{\pi_{n+1}} ( \Delta )$ $< \Prob_{\pi_{n}} ( \Delta ) = \delta$,
contradicting the minimality of $\delta$.
Therefore $\delta=0$ and $\mu=\mu_d$, as claimed.

This completes the proof of Theorem 3.

We remark that the preceding demonstration 
is very similar to the one expounded 
in sections 2 and 3 of chapter VIII of \cite{liggett}
in the case of continuous time.
However, we have made the proof slightly easier
by invoking compactness in order to choose
a discrepancy-minimizing $\pi$ at the outset
(since we then need only show that
we can decrease the discrepancy further,
rather than show that we can reduce it all the way to zero).
Also, we have coupled $\mu$ directly with $\mu_d$,
rather than with $\mu_{d'}$
for $d'$ close to $d$.

\section{Monotonicity Properties of Particle Density}
\label{sec:monotonicity}

Before we dive into the technicalities that will ultimately afford us
a rigorous proof of Theorems 1, 2, and 3,
we offer a heuristic reason for the circularity of the temperate zone,
in terms of the behavior of the particle process.

We will assume heuristically that short excerpts from the lattice-path
exhibit statistics that are nearly Markovian
and are governed by some particular measure $\mu_d$,
where $d$ is not constant but is a slowly varying function of position.

Return briefly to the first-quadrant Young-diagram picture.
If one examines a piece of the lattice-path
that is approximately governed by some particular $\mu_d$,
one sees that the path approximates a line of slope
$m = -\frac{p_1}{p_0} = -\frac{d}{1-d} < 0$.
Write this line as $y=mx+b$.
If every point on this line moved 1 unit to the right and 1 unit upward,
the equation of the line would become $y-1=m(x-1)+b$
or $y=mx+b+(1-m)$;
that is, the line would move upward by $1-m$.
Similarly, if every point on the lattice-path
moved 1 unit to the right and 1 unit upward,
the straight line that locally approximates the path
would appear to slide $1-m$ units upward.
However, when we perform the stochastic lattice-path update,
only some of the points get moved
1 unit to the right and 1 unit upward,
namely, a proportion of $\frac12 p_1 q_{10}$
(where $p_1 q_{10}$ is the frequency of down-then-right bends
and $\frac12$ is the probability that a given down-then-right bend
will become a right-then-down bend, or equivalently,
that the halfway point along that jag
will move 1 unit to the right and 1 unit upward).
Hence the lattice-path is expected to move 
a mean upward distance of only
$\frac12 p_1 q_{10} \cdot (1-m)$.
But $p_1 q_{10} = p_0 q_{01}$
and $1-m=1+\frac{p_1}{p_0}=\frac{1}{p_0}$,
so the expected vertical displacement is
$\frac12 p_0 q_{01} \cdot \frac{1}{p_0} = \frac12 q_{01}$.
We wish to express this in terms of $m$,
so first we rewrite it in terms of $d$.
Recall $q(0,1) = \frac{1-s}{1-d}$
with $s = \sqrt{d^2+(1-d)^2}$.
This yields 
\begin{align*}
q(0,1)
&= \frac{1}{1-d} - \frac{s}{1-d} \\
&= \frac{1}{1-d}
   - \sqrt{1 + \bigl(\tfrac{d}{1-d}\bigr)^2};
\end{align*}
performing the substitution $d = -\frac{m}{1-m}$,
we obtain 
\begin{equation*}
\frac12 q(0,1) = \frac{1-m-\sqrt{1+m^2}}{2}.
\end{equation*}

Now let us assume that when $n$ is large,
the lattice-path is an approximation to some continuous curve
$\frac{Y}{n} = \phi(\frac{X}{n})$,
where $\phi$ satisfies boundary conditions 
$\phi(0)=\frac12$ and $\phi(\frac12)=0$.
Note that under this assumption,
the lattice-path in the vicinity of the point $(X,Y)$
(with $\frac{Y}{n} = \phi(\frac{X}{n})$)
should drift upward by approximately
\begin{align*}
(n+1)\phi ({\scriptstyle \frac{X}{n+1}}) - n\phi ({\scriptstyle \frac{X}{n}})
&\approx (n+1) \bigl( \phi ({\scriptstyle \frac{X}{n}}) 
   -  ({\scriptstyle \frac{X}{n} - \frac{X}{n+1}}) 
      \phi'({\scriptstyle \frac{X}{n}}) \bigr) 
   - n\phi({\scriptstyle \frac{X}{n}}) \\
&= \phi({\scriptstyle \frac{X}{n}}) 
     - {\scriptstyle \frac{X}{n}} \phi'({\scriptstyle \frac{X}{n}}) \\
&= {\scriptstyle \frac{Y}{n}} 
     - {\scriptstyle \frac{X}{n}} \phi'({\scriptstyle \frac{X}{n}}).
\end{align*}

Putting 
$x=\frac{X}{n}$,
$y=\phi(x)=\frac{y}{n}$, and
$\phi'(\frac{X}{n}) = s = \frac{dy}{dx}$,
and equating the two drift rates,
we get the differential equation
\begin{equation*}
y-x\frac{dy}{dx} = \frac{1-\frac{dy}{dx}-\sqrt{1+(\frac{dy}{dx})^2}}{2}
\end{equation*}
for $0 \leq x \leq \frac12$.
It is easy to check that 
$y=\frac12-\sqrt{x-x^2}$ is a solution;
but this is just the lower-left quarter-circular arc
of the circle
$(x-\frac12)^2+(y-\frac12)^2=\frac14$.

It would be satisfying to make the preceding heuristic argument 
rigorous using general hydrodynamic-limit methods for 
exclusion-type processes, but we did not pursue that approach here; 
instead, we followed the technique developed by Rost \cite{rost} for
the continuous-time model and adapted it to our discrete-time setting.
In the present section we are mainly interested in
the limiting behavior of the $X$-process as one moves out from
the origin along a space-time line of slope $u$, i.e., the pictures
seen near location $un$ at time $n$.  To study this, we introduce
a function $h(\cdot)$ with the property that the number of particles 
to the right of position $k$ at time $n$ is roughly $n\:h(k/n)$.  
We then show that $h$ exists almost everywhere and is a convex 
(hence almost everywhere differentiable) function of $u$ 
(Proposition~\ref{prop2}).  We further show that the derivative exists, 
and is essentially the negative of the probability of finding a 
``1'' at locations near $un$ at time $n$, where $n$ is large 
(Proposition~\ref{prop3}).  We then show that at a fixed time, the local
statistics change monotonely as one moves to the right, in the sense
that any given pattern of 1's is less likely to appear the further
one goes to the right (Proposition~\ref{prop4}).  This result is used to show
that the limiting behavior mentioned above is almost everywhere a 
convex combination of the update- and translation-invariant measures 
studied in the previous section (Proposition~\ref{prop5}).  Finally, we show
that these local statistics depend in a continuous way on the speed
$u$ of the ``observer'' (Proposition~\ref{prop6}).

In the section following this one, we obtain a formula for $h$
by means of two separate arguments, which give upper and lower
bounds for $h$ that turn out to coincide.
We obtain a lower bound on $h$ by slowing down the lead particle,
which clearly cannot increase the function $h$.
We obtain an upper bound by looking at how fast the
``1''s can move when they have a given density; the essential idea is that
at location $un$ and time $n$, an upper bound on 
this ``average velocity'' can be
calculated and eventually leads to an upper bound on $h(u)$.
	
We have made no attempt to hide
our indebtedness to Rost's work;
whole paragraphs have been lifted from his article
with only minor modifications.
We have done this 
because we do not see room
for many improvements in Rost's exposition,
and because we wanted our presentation to be self-contained.

As in \cite{rost}, we will let $X$ be our particle process
with state space $\cX = \{0,1\}^\Z$,
except that the updates only occur at discrete moments
indexed by $\N$.
$X(k,n)$ is the state of position $k$ at time $n$,
with $k$ in $\Z$ and $n$ in $\N$.
The initial state ($n=0$) is given by 
$X(k,0) = x_k^*$, where
\[
x_k^* = \left\{ \begin{array}{ll}
	1 & \text{if $k \leq 0$,} \\
	0 & \text{if $k > 0$.}
	\end{array} \right.
\]
The particle initially at location 0 will be called the {\bf lead particle}.
The order on points $x,y \in \cX$ defined by
\[
x \leq y \ \ \ \text{if and only if} 
\ \ \ x_i \leq y_i \ \text{for all $i \in \Z$}
\]
induces a {\bf stochastic order} on the set of probability 
measures on $\cX$, namely
$\mu \leq \nu$ if and only if $\pi(\{(x,y): \ x \leq y \}) = 1$
for some joining $\pi$ of $\mu$ and $\nu$,
where a joining of two probability measures $\mu$, $\nu$ on $\cX$
is a probability measure on $\cX^{(1)} \times \cX^{(2)}$
with respective marginals $\mu$ and $\nu$.
A $k$-point correlation of $\mu$ (with $k \geq 1$)
is a quantity of the form
$\mu(\{x:\ x_{i_1} = x_{i_2} = \dots = x_{i_k} = 1 \} )$,
with $i_1 < i_2 < \dots < i_k$;
a measure $\mu$ on $\cX = \{0,1\}^\Z$ is determined by its
correlations~\cite{billingsley}.

Note that, in the specified initial state $x^*$
and in all states accessible from $x^*$
in finite time,
there is a rightmost 1.
We may therefore define the process
\[
S(k,n) = \sum_{i>k} X(i,n)
\]
with state space $\cS$ equal to
the set of all (weakly) decreasing sequences of non-negative integers,
and we may further define a component-wise order on $\cS$
and a stochastic order on the set of probability measures on $\cS$
just as we did for $\cX$.

We let $\cL(Y)$ denote the probability law
governing a random variable $Y$,
we let $\cE(Y)$ denote the expected value
of $Y$, and
we let $*$ denote convolution of probability measures on $\N$.

\begin{proposition}
\label{prop1}
For all non-negative integers $m,n$ and integers $k,l$, 
one has 
\begin{equation*}
\cL(S(k,m)) * \cL(S(l,n)) \geq \cL(S(k+l,m+n)).
\end{equation*}
\end{proposition}

\begin{proof}
A simple coupling argument shows
that the evolution-rule of the $S$-process
preserves stochastic order of measures on its state space $\cS$.
(In the interpretation of particle-configurations
as lattice-paths,
one state of the $S$-process is dominated by another
if the lattice-path associated with the former
never crosses above or to the right of
the lattice-path associated with the latter.)
Accordingly, we define a process $\tilde{S}$
that is equal to $S$ up until time $m$
but at time $m$ is replaced by
\[
\tilde{S}(j,m) = \left\{ \begin{array}{ll}
	S(k,m) & \text{for $j \geq k$ and} \\
	S(k,m)+(k-j) & \text{for $j < k$;} 
\end{array} \right.
\]
that is, all the particles that are to the left of position $k+1$
at time $m$
simultaneously move as far to the right as possible,
so that site $k$ and all sites to its left become occupied.
After time $m$, $\tilde{S}$ again evolves in accordance
with the dynamics of the $S$-process.
One sees that $\cL(\tilde{S}(\cdot,m+n)) \geq \cL(S(\cdot,m+n))$.
But, conditioned on $S(k,m)$,
the law of $\tilde{S}(k+l,m+n)-S(k,m)$ ($l \in \Z$, $n \geq 0$)
is identical to that of $S(l,n)$ ($l \in \Z$, $n \geq 0$)
and independent of $S(k,m)$.
Therefore $\cL(\tilde{S}(k+l,m+n))=\cL(S(k,m))*\cL(S(l,n))$.
\end{proof}

\begin{proposition}
\label{prop2}
For all $u \in \R$,
the random variables $\frac1n S(\lfloor un \rfloor,n)$
converge almost surely and in $L^1$ to
a constant $h(u)$ as $n \rightarrow \infty$.
The function $h$ is decreasing and convex;
$h(u)=0$ for $u > \frac12$ and
$h(u)=-u$ for $u < -\frac12$.
\end{proposition}

\begin{proof}
Proposition~\ref{prop1}, coupled with the fact that
$\lfloor um \rfloor + \lfloor un \rfloor \leq \lfloor u(m+n) \rfloor$,
gives us
\[ \cL(S(\lfloor um \rfloor,m)) * \cL(S(\lfloor un \rfloor,n))
\ge \cL(S(\lfloor u(m+n) \rfloor,m+n)).  \]
The convergence statements in the Proposition
follow from the Kesten-Hammersley theorem \cite{smythe}.
To prove convexity of $h$, we deduce from Proposition~\ref{prop1}
that for $\alpha,\beta > 0$ with $\alpha+\beta=1$,
\[
\cE S(\lfloor \alpha u n \rfloor, \alpha n) +
\cE S(\lfloor \beta  v n \rfloor, \beta  n) \geq
\cE S(\lfloor ( \alpha u + \beta v ) n \rfloor, n).
\]
Dividing both sides by $n$ gives
\[
\alpha h(u) + \beta h(v) \geq h ( \alpha u + \beta v ).
\]
If $u > \frac12$, we have $h(u)=0$,
as the lead particle moves to the right with mean speed $\frac12$.
Similarly, the leftmost vacancy moves to the left with mean speed $-\frac12$,
whence $h(u)=-u$ for $u < -\frac12$.
\end{proof}

\begin{proposition}
\label{prop3}
If $h$ is differentiable at $u$, then
$\cE X(k,n) \rightarrow -h'(u)$
whenever $n \rightarrow \infty$ with $k/n \rightarrow u$.
\end{proposition}

\begin{proof}
We consider the functions $h_n$, defined by
\begin{equation*}
h_n (v) = \int_v^\infty \cE X ( \lfloor wn \rfloor, n ) \ dw;
\end{equation*}
note that $h_n(v)$ is roughly equal to
$\frac{1}{n} \sum_{k=\lfloor vn \rfloor}^\infty \cE X(k,n)$.

Since for each fixed $n$
the function $k \mapsto \cE X(k,n)$ is decreasing in $k$,
the function $v \mapsto \cE X(\lfloor vn \rfloor,n)$ 
is decreasing in $v$, implying that $h_n(v)$,
being the integral of a decreasing function, is convex.
The functions $h_n(v)$ tend to $h(\cdot)$ as $n$ gets large, where $h(v) = 
\lim_{n \rightarrow \infty} \frac{1}{n} \cE S ( \lfloor v n \rfloor, n)$.
Hence, by an elementary lemma on convex functions of a real argument,
the desired result holds.
Indeed, more is true;
as long as $h'(u)$ exists,
$h'_n (v) \rightarrow h'(u)$ as $v \rightarrow u$,
where $h'_n(v)$ may be either the right or left derivative of $h_n$ at $v$.
\end{proof}

Define $f(u) = -h'(u+0)$, where $h'(u+0)$ signifies
the right derivative of $h$ at $u$.
$f(u)$ is our candidate for the density of the particle process
at location $\lfloor un \rfloor$ at time $n$,
when $n$ is large.

Let $\mu(k,n) = \cL(X(k+l,n), l \in \Z)$.
That is, $\mu(k,n)$ is the probability law governing $X(\cdot,\cdot)$,
shifted $k$ positions spatially and $n$ steps into the future.

\begin{proposition}
\label{prop4}
For all $k \in \Z$ and $m,n \in \N$
we have
\begin{align}
\mu(k,n) &\geq \mu(k+1,n), 
	\label{twelve} \\
\mu(k,n+m) &\leq \sum_l \beta(m,l) \mu(k-l,n), \ \text{and} 
	\label{thirteena} \\
\mu(k,n+m) &\geq \sum_l \beta(m,l) \mu(k+l,n), 
	\label{thirteenb}
\end{align}
where $\beta(m,\cdot)$ is the binomial distribution
with mean $\frac m 2$ and variance $\frac m 4$.
\end{proposition}

\begin{proof}
$\mu(k+1,n)$ is the law of $X(k+l,n)$, $l \in \Z$
under the initial condition $x'$,
where $x'_i$ is 1 or 0 according to whether 
$i \leq -1$ or $i > -1$;
since $x' \leq x^*$,
\eqref{twelve} follows from the monotonicity of the dynamics.

Monotonicity also plays a role in the proof of \eqref{thirteena}.
The position of the lead particle at time $m$
is binomially distributed with mean $\frac{m}{2}$;
conditioned on that position,
one compares the original process with the process
in which all sites behind the first particle are occupied at time $m$,
and which evolves according to the usual dynamics after time $m$.
\eqref{thirteenb} follows from \eqref{thirteena},
by symmetry between migration of particles and migration of vacancies.
\end{proof}

\begin{proposition}
\label{prop5}
If $h$ is differentiable at $v$,
any weak limit $\mu^*$ of the measures $\mu(\lfloor un \rfloor , n)$, 
$n \rightarrow \infty$, is of the form
\begin{equation*}
\mu^* = \int_0^1 \mu_t \, \rho(dt)
\end{equation*}
with $\rho$ some probability distribution on $[0,1]$,
and with $\mu_t$ defined as in section~\ref{sec:reduction}~(there we called it $\mu_d$).
\end{proposition}

\begin{proof}
$\mu^*$ is stochastically larger than its image under the shift,
by Proposition~\ref{prop4}.
Since 
$\frac{\lfloor un \rfloor + k}{n} \rightarrow u$ as $n \rightarrow \infty$,
for every fixed $k$,
both $\mu^*$ and its image under the shift
have the same one-point correlations, namely $f(u)$
(Proposition~\ref{prop3}).
It follows that $\mu^*$ and its image under the shift are identical.
But shift-invariance of $\mu^*$
implies also its invariance under the action
of the time-evolution semigroup, if one uses 
the second and third inequalities of Proposition~\ref{prop4}.
Using the main result of the previous section,
we arrive at the desired result.
\end{proof}

If $F$ is a finite subset of $\Z$,
let $\rho(k,F;n)$ be the probability that
$X(k+i,n)=1$ for all $i \in F$.

\begin{proposition}
\label{prop6}
Assume $h$ is differentiable at $\ou$.
For any finite set $F$ of cardinality $n$ and any $\epsilon > 0$
there exists a $\delta > 0$ and $n_0 \in \N$ such that
\begin{equation}
\label{sixteen}
|\rho(\lfloor un \rfloor,F;n)-\rho(\lfloor \ou n \rfloor,F;n)| 
\leq \epsilon
\end{equation}
for $|u - \ou | \leq \delta$ and $n \geq n_0$, and
\begin{equation}
\label{seventeen}
|\rho(\lfloor \ou n \rfloor,F;n+m)-
\rho(\lfloor \ou n \rfloor,F;n)| 
\leq \epsilon
\end{equation}
for $0 \leq m \leq \delta n$ and $n \geq n_0$.
\end{proposition}

\begin{proof}
First let us prove \eqref{sixteen}.
By symmetry, it suffices to handle the case $u>\ou$.

Since $\mu(k,n)$ is stochastically decreasing in $k$
we need an upper estimate only for
\begin{equation}
\label{nolabel}
\rho(\lfloor \ou n \rfloor , F; n)
 - \rho(\lfloor u n \rfloor , F ; n) .
\end{equation}
But the definition of stochastic order (via coupling) 
gives us such an estimate:
\begin{equation}
\sum_{i \in F} \bigl( \rho( \lfloor \ou \rfloor , \{i\} ; n)
                    - \rho( \lfloor u \rfloor , \{i\} ; n) \bigr)
= \sum_{i \in F} 
\cE(X(\lfloor \ou n \rfloor + i , n) - X(\lfloor u n \rfloor + i , n)).
\end{equation}
Take $\delta>0$ such that $h'(\ou+\delta)$ exists
and satisfies
$-h'(\ou+\delta)=f(\ou+\delta) \geq f(\ou) - \frac{\epsilon}{2n}$
(recall that $f$ is assumed continuous at $\ou$).
Proposition~\ref{prop3} gives
\begin{equation}
\liminf_{n \rightarrow \infty} \cE(X(\lfloor (\ou+\delta) n \rfloor + i, n)
\geq f(\ou) - \frac{\epsilon}{2n}.
\end{equation}
Hence we have, uniformly for $u \leq \ou + \delta$,
\begin{equation}
\limsup_{n \rightarrow \infty} \sum_{i \in F}
\bigl( \rho(\lfloor \ou n \rfloor , \{i\}; n)
- \rho(\lfloor u n \rfloor , \{i\}; n) \bigr)
\leq \epsilon/2.
\end{equation}
This implies that \eqref{nolabel}
becomes eventually smaller than $\epsilon$,
uniformly in $u \leq \ou + \delta$.

This completes the proof of \eqref{sixteen}.
The estimate \eqref{seventeen}
follows from \eqref{sixteen}
in combination with Proposition~\ref{prop4},
inequalities \eqref{thirteena} and \eqref{thirteenb},
and the fact that $\beta(s,\cdot)$
is essentially carried by a set of the form
$\{l: l \leq cm\}$ in the limit $m \rightarrow \infty$.
\end{proof}

\section{Identification of the Density Profile}
\label{sec:profile}

\subsection{Lower bound}
\label{subsec:lower}

Let $Z(i,n)$ be the position at time $n$
of the particle originally located at position $-i$,
and for $k \geq 1$
let $Y(k,n) = Z(k-1,n) - Z(k,n)$, 
the distance between the $k-1$st and $k$th particles from the right
at time $n$.
Unfortunately, the expected value of $Y(k,n)$
goes to infinity as $n$ gets large, with $k$ fixed,
so that the law of large numbers cannot be applied
as directly as we might like.
Rost was able to remedy this problem in a clever way,
as we are about to see.

\begin{proposition}
\label{prop7}
$h(u) \geq \frac{1-u}{2} - \frac12 \sqrt{\frac12 - u^2}$
for $|u| \leq \frac12$.
\end{proposition}

\begin{proof}
We modify the dynamics of the system as follows:
when deciding whether to advance the lead particle or not,
we use a biased coin,
so that it advances with probability $\frac b 2$ with $b \leq 1$;
the other particles advance with probability $\frac 1 2$ as before.
Expectations with respect to this process will be denoted by $\cE^b$,
and its probability law will be denoted by $P^b$.

These dynamics, in terms of the $Y$-process,
may be described as follows:
The state space is the set of all sequences of positive integers
$(y_1,y_2,\dots)$.
Given such a sequence,
let $\alpha_1$ be a Bernoulli random variable with expected value $\frac b 2$,
and for $i \geq 2$ let $\alpha_i$ be a Bernoulli random variable
of expected value $\frac 1 2$
unless $y_{i-1} = 1$,
in which case let $\alpha_i$ be the constant 0.
(We can think of $\alpha_i$ as the indicator function of
the event in which the $i$th particle moves to the right.)
Then the $Y$-process, when in state $(y_1,y_2,\dots)$,
jumps to $(y_1+\alpha_1-\alpha_2,y_2+\alpha_2-\alpha_3,\dots)$.
One can check two statements about this process:
first, it preserves stochastic order;
second, an invariant measure is $\gamma^b$,
defined by the properties
that all its coordinates are 
independent and identically distributed
and that
\begin{equation*}
\gamma^b(y_i>m)= \left\{ \begin{array}{ll}
	1 & \text{if $m=0$}, \\
	b \bigl( \frac{b}{2-b} \bigr) ^{m-1} & \text{if $m>0$}.
	\end{array} \right.
\end{equation*}

If we compare the modified $Y$-process
with initial condition $y_1=y_2=\dots=1$
and the modified $Y$-process
with initial condition given by the stationary measure $\gamma^b$,
we find that the law of the first process
is stochastically smaller than
the law of the second process at time $n=0$
and hence for all times.
We thus get, for all $k \geq 1$ and $n \geq 0$,
\begin{equation*}
\cE^b \bigl( \sum_{j=1}^k Y(j,n) \bigr)
\leq k \sum_{m \geq 0} \gamma^b (y_i > m) = k \: \frac{2-b^2}{2-2b}.
\end{equation*}
If we choose $k= \lfloor an \rfloor$, with $0 \leq a \leq 1$ fixed, 
and let $n \rightarrow \infty$, we get
\begin{equation*}
\limsup_{n \rightarrow \infty} \ \frac{1}{n} 
\cE^b \Bigl( \sum_{j \leq an} Y(j,n) \Bigr) \leq a \: \frac{2-b^2}{2-2b}.
\end{equation*}
In fact, using the weak law of large numbers for $\gamma^b$
and stochastic domination,
we find that for any $\epsilon > 0$,
\begin{equation*}
P^b \Bigl[ \frac1n \sum_{j \leq an} Y(j,n) > 
a \: \frac{2-b^2}{2-2b} + \epsilon \Bigr] \rightarrow 0
\end{equation*}
as $n \rightarrow \infty$.
For now, hold both $\epsilon$ and $a$ fixed.
We know that $Z(0,n)$ is binomially distributed with mean $\frac b 2$,
so the law of large numbers, in combination with the preceding result
and the fact that $Z(\lfloor an \rfloor,n) = Z(0,n)-\sum_{j \leq an} Y(j,n)$,
gives
\begin{equation*}
P^b \Bigl[ \frac1n Z(\lfloor an \rfloor,n) < 
\frac b 2 - a \: \frac{2-b^2}{2-2b} - \epsilon \Bigr] \rightarrow 0.
\end{equation*}

The $Z$-process gets stochastically larger if $b$ is replaced by 1;
hence the preceding proposition remains true
if we replace $P^b$ by $P$,
the probability law governing the original dynamics of the $Z$-process.
This holds for every $b$.
In particular, setting $b=1-\sqrt{\frac{a}{1-a}}$, so that
\begin{equation*}
\frac{b}{2} - a \: \frac{2-b^2}{2-2b} = \frac12 - a - \sqrt{a(1-a)}
\end{equation*}
(the maximum value of $\frac b 2 - a \frac{2-b^2}{2-b}$
as $b$ ranges over $[0,1]$),
we get
\begin{equation*}
P \Bigl[ \frac1n Z(\lfloor an \rfloor, n) 
< \frac12 - a - \sqrt{a(1-a)} - \epsilon \Bigr]
\rightarrow 0.
\end{equation*}
Expressing this in terms of the $S$-process one obtains
\begin{equation*}
P \Bigg[ \frac1n 
S\bigg(\Big\lfloor \big( \frac12-a-\sqrt{a(1-a)}-\epsilon \big) n 
\Big\rfloor , \ n \bigg) 
> a \Bigg] \rightarrow 1.
\end{equation*}
Since this is true for all $\epsilon > 0$,
\begin{equation*}
h\Big(\frac12 - a - \sqrt{a(1-a)}\Big) \geq a.
\end{equation*}
Setting $u=\frac12 - a - \sqrt{a(1-a)}$,
we conclude that
$h(u) \geq \frac{1-u}{2} - \frac12 \sqrt{\frac12 - u^2}$
for $|u| \leq \frac12$, as claimed.
\end{proof}

\subsection{Upper bound}
\label{subsec:upper}

To complement Proposition~\ref{prop7}, we have

\begin{proposition}
\label{prop8}
$h(u) \leq \frac{1-u}{2} - \frac12 \sqrt{\frac12 - u^2}$
for $|u| \leq \frac12$.
\end{proposition}

\begin{proof}
First assume $u>0$, and put $w=1/u$.
Assume that $u$ is irrational and that $h'(u)$ exists
(as must be the case for a dense set of $u$'s in $(0,\frac12)$).
We compute the expected value of $S(\lfloor un \rfloor,n)$, which is
the expected value of the number of particles
that have passed an observer traveling at speed $u$,
minus the number of particles that the observer has passed:
\begin{align*}
\cE S(k, \lfloor kw \rfloor)
&= \tfrac12 \sum_{i=0}^{\lfloor kw \rfloor}
   P\bigl[
     X(\lfloor iu \rfloor, i) = 1,\;
     X(\lceil iu \rceil, i) = 0
   \bigr] \\
&\quad - \sum_{l=1}^k \cE X(l, lw).
\end{align*}
Multiplying both sides by $u/k \approx 1/\lfloor kw \rfloor$ 
and taking the limit, one gets
\begin{align*}
h(u)
&= \lim_{k \rightarrow \infty} \frac{\cE S(k,\lfloor kw \rfloor)}
{\lfloor kw \rfloor} \\
&= \frac12 \Big( \lim_{k \rightarrow \infty}
\frac{1}{\lfloor kw \rfloor} \sum_{i=0}^{\lfloor kw \rfloor} P [
X(\lfloor iu \rfloor , i) = 1, \ X(\lceil iu \rceil , i) = 0 ] \Big) \\
& \quad - \ \frac1w \Big( \lim_{k \rightarrow \infty}
\frac1k \sum_{l=1}^k \cE X(l,lw) \Big) \\
&= \frac12 \lim_{N \rightarrow \infty} \frac1N \sum_{i=0}^{N-1}
\mu(\lfloor ui \rfloor,i) (x_0=1,\ x_1=0) \ -uf(u).
\end{align*}
The same argument shows that this formula holds
for negative $u$, too.
Thus this relation holds for a dense set of $u$'s
in $[-\frac12,\frac12]$.

Now, by Proposition~\ref{prop5}, any limit of
\begin{equation*}
\mu(\lfloor un \rfloor,n)(x_0=1,\ x_1=0)
\end{equation*}
as $n \rightarrow \infty$ is of the form
\begin{equation*}
\int (1-\sqrt{a^2+(1-a)^2}) \ \rho(da) \ \ \text{with} 
\ \ \int a \ \rho(da) = f(u)
\end{equation*}
(see the formula for the probability of the event $x_0=1$, $x_1=0$
calculated at the end of subsection~\ref{subsec:markov}).
Hence by Jensen's inequality we get
\begin{equation*}
\limsup_{n \rightarrow \infty} \mu(\lfloor un \rfloor,n)(x_0=1,\ x_1=0)
\leq 1-\sqrt{(f(u))^2+(1-f(u))^2}
\end{equation*}
so that
\begin{equation*}
{\scriptstyle \frac12} \lim_{N \rightarrow \infty} \sum_{i=0}^{N-1} 
\mu(\lfloor ui \rfloor,i)(x_0=1,\ x_1=0)
\leq {\scriptstyle \frac12}-{\scriptstyle \frac12} \sqrt{(f(u))^2+(1-f(u))^2}.
\end{equation*}
Hence
\begin{equation*}
h(u) \leq {\scriptstyle \frac12}-
{\scriptstyle \frac12} \sqrt{(f(u))^2+(1-f(u))^2} - uf(u).
\end{equation*}
Since $0 \leq f(u) \leq 1$, we have
\begin{equation*}
h(u) \leq \sup_{0 \leq b \leq 1} \Big\{ {\scriptstyle \frac12} - 
{\scriptstyle \frac12} \sqrt{b^2+(1-b)^2} - ub \Big\}.
\end{equation*}
This is maximized at $b = \frac12 - \frac{u}{\sqrt{2-4u^2}}$ (=$f(u)$).
Substituting, we get
$$h(u) \leq \frac12 - \frac12 \sqrt{\frac12 - u^2} - \frac{u}2.$$
\end{proof}

\subsection{Conclusion of Proof}
\label{subsec:finish}

Propositions 2, 3, 7, and 8
combine to yield the following density profile and law of large numbers:
For any $u \in \R$, 
$\cE X(k,n)$ tends towards $f(u)$ 
as $n$ goes to infinity with $k/n$ tending towards $u$.
The function $f$ is given by
\[
f(u) = \left\{ \begin{array}{ll}
	1 & 
		\text{for $u<-\frac12$,} \\
	\frac12 - \frac{u}{\sqrt{2-4u^2}} & 
		\text{for $-\frac12 \leq u < \frac12$, and} \\
	0 & \text{for $u>\frac12$}. \end{array} \right.
\]
The quantities $\frac1n \sum_{un<k<vn} X(k,n)$ converge almost surely
to the constant value $$\int_u^v f(w) \ dw$$ for $u < v$.

We can now unwind our results to obtain a proof of
the Arctic Circle Theorem.
Under a change in coordinates,
the preceding result says that if
we evolve an infinite lattice-path in the first quadrant
in the fashion described near the end of section~\ref{sec:reduction},
the lattice-path at time $n$
will almost surely attach to the $x$- and $y$-axes
at points $(\frac{n}2 + o(n),\:0)$
and $(0,\:\frac{n}2 + o(n))$,
and for all $0 < \theta < \frac\pi2$,
the lattice-path at time $n$ will almost surely
cross the line $y/x = \tan \theta$ at a point
$(\frac{n}2 \cos \theta + o(n),
\frac{n}2 \sin \theta + o(n))$.
For each fixed $\epsilon>0$,
we can find angles
$0 < \theta_1 < \theta_2 < \dots < \theta_m < \frac\pi2$
and a number $\delta>0$
so that any curve that starts at $(0,\frac{n}2 \pm \delta n)$,
ends at $(\frac{n}2 \pm \delta n,0)$,
only moves rightward and downward,
and meets each line $y/x = \tan \theta_i$
within distance $\delta$ of
$(\frac{n}2 \cos \theta_i, \frac{n}2 \sin \theta_i)$
(for all $1 \leq i \leq m$)
must necessarily stay within distance $\epsilon$
of the quarter-circle $x^2 + y^2 = \frac{n^2}{4}$,
$x,y \geq 0$.
This concludes the proof of the Arctic Circle Theorem.

It is also worth pointing out that the measures $\mu_d$
are in fact ``attractors'' 
under the time-evolution map $F:\cM(\cX) \rightarrow \cM(\cX)$
if one restricts to ergodic translation-invariant measures
on $\{0,1\}^{\Z}$
(that is, those translation-invariant measures that cannot be decomposed
as convex combinations of other translation-invariant measures).
For, let $\mu'$ be any such measure,
with $\mu'(\{x:x_0=1\})=d$,
and take a measure $\pi_0$ on $\cX^{(1)} \times \cX^{(2)}$
with marginals $\mu_d$ and $\mu'$.
Let $\pi_i = \overline{F}^i(\pi_0)$
under the coupling-dynamics 
$\overline{F}: \cM(\cX^{(1)} \times \cX^{(2)}) 
\rightarrow \cM(\cX^{(1)} \times \cX)$.
We have shown that $\pi_i$ converges almost surely to a diagonal measure
with projection $\mu_d$ on each component;
hence the $\pi_i$'s converge in distribution, and in particular,
$F^i(\mu')$ (the second marginal of $\pi_i$) converges to $\mu_d$.

\section{Introducing Bias}
\label{sec:bias}

Domino tilings have been studied 
by researchers in statistical mechanics
in another guise,
namely, dimer-patterns on a square grid;
see for instance \cite{kasteleyn}.
Suppose that, as in Kasteleyn's paper,
we introduce an energy function that discriminates
between horizontal and vertical tiles.
It is shown in \cite{elkies} that
the Aztec diamond of order $n$ has $n(n+1)/2 \choose k$ tilings
with $2k$ horizontal tiles and $n(n+1)-2k$ vertical tiles
(and no tilings in which the number of horizontal tiles
or vertical tiles is odd).
Thus, if we fix $0 \leq p \leq 1$ and 
assign measure $p^k (1-p)^{n(n+1)/2-k}$
to each tiling of the Aztec diamond of order $n$
with $2k$ horizontal tiles and $n(n+1)-2k$ vertical tiles,
we obtain a probability measure on the set of tilings
that is a Gibbs state relative to an energy function
that assigns energy $-\frac12 \log p$ to horizontal dominoes
and $-\frac12 \log (1-p)$ to vertical dominoes.

It is not hard to show
that random domino tilings of Aztec diamonds
under this energy function
may be iteratively generated by ``biased shuffling,''
in which an empty 2-by-2 block is filled
with two horizontal dominoes with probability $p$
and two vertical dominoes with probability $1-p$.
The temperate zone can be defined as before,
and it was shown in \cite{CEP} that its boundary is 
the ellipse $\frac{x^2}{p}+\frac{y^2}{1-p}=1$.

We will prove part of this assertion rigorously,
using the methods introduced earlier.
Each of the four growing ``frozen'' regions
can be associated with a Ferrers-diagram (or Young-diagram) growth process.
For the north and south frozen regions,
the probability of growth is $p$;
for the east and west, it is $1-p$.
In all cases, the growth process may in turn be replaced
by an asymmetric exclusion process,
in which the probability of a particle moving to its right
(assuming that there is a vacancy there)
is equal to the growth-rate for the growth-process.

As in the unbiased case,
there exists (for each value of the bias $p$)
a one-parameter family of 
extremal elements in the set of 
dynamically stationary, shift-invariant measures,
where the parameter corresponds to density;
more specifically, we get Markov measures on $\cX = \{0,1\}^\Z$
with transition probabilities
$q_{01} = \frac{1-\sqrt{1-4pd(1-d)}}{2p(1-d)}$,
$q_{10} = \frac{1-\sqrt{1-4pd(1-d)}}{2pd}$,
$q_{00}=1-q_{01}$,
$q_{11}=1-q_{10}$.
Here $p$ is the growth rate (or drift rate)
and $d$ is the density of 1's.
The heuristic method given at the beginning of section~\ref{sec:monotonicity}~can
be used to lend support to our ellipse conjecture.
If we knew that these were the only 
dynamically stationary, shift-invariant measures,
we would be able to deduce an
``arctic ellipse theorem''.
Unfortunately, when $p$ exceeds one-half,
the coupling method of section~\ref{sec:steady}~does not work.
Specifically, it is no longer possible
to devise a coupling under which
the number of mismatches
is guaranteed to weakly decrease over time.
One can see this by noting that if $p$ is close to 1,
then the block
\[
\begin{array}{lll}
        1 & 0 & 0\\
        1 & 1 & 0
\end{array}
\]
is very likely to give rise to the block
\[
\begin{array}{lll}
        0 & 1 & 0\\
        1 & 0 & 1
\end{array}
\]
at the next time-step, 
so that the number of mismatches within the block 
will increase from 1 to 3.
If the initial state of the joint process is
\[
\begin{array}{llllllllllllll}
        \dots & 1 & 0 & 0 & 1 & 1 & 0 & 1 & 0 & 0 & 1 & 1 & 0 & \dots \\
        \dots & 1 & 1 & 0 & 1 & 0 & 0 & 1 & 1 & 0 & 1 & 0 & 0 & \dots
\end{array}
\]
then the density of mismatches will increase from $1/3$ to nearly $1$.

It is nevertheless true that when $p$ is less than one-half,
the proof of Theorem 3 given above
carries over in a straightforward way,
showing that the Markov measures introduced in the preceding paragraph
(with $p$ fixed and $d$ varying)
are the only stationary translation-invariant probability measures.
Hence, the north and south frozen regions
are indeed bounded by the prescribed arcs of ellipses.
In the case where $p$ is greater than one-half,
the heuristic method described at the beginning of 
section~\ref{sec:monotonicity}~leads one to a differential equation
that is satisfied by the ellipse, and later work (see \cite{CEP})
confirmed that this is in fact the correct answer.

It is interesting to note that if one takes the limiting behavior of the
biased exclusion process in discrete time as $p \rightarrow 0$,
the elliptical arc tends in shape to a parabolic arc
(once a renormalization is made to compensate for
the fact that the arc is getting smaller).
This is reassuring, since the discrete-time exclusion process
should ``converge'' to the continuous-time exclusion process
as $p \rightarrow 0$ when time is re-scaled,
and since (as was remarked earlier)
the continuous-time process
was shown by Rost to yield a parabolic arc as its asymptotic profile.

\section{Conclusion}
\label{sec:conclusion}

We have shown that a random domino-tiling of a large Aztec diamond
(unlike a random domino-tiling of a large square)
is likely to have large-scale structure,
and we have given precise information about this structure.
Outside of a certain critical circle
(which we call the arctic circle),
a random tiling is likely to exhibit four different sorts of local behavior,
associated with four particular tilings of the plane.

A natural next step is describing
the behavior of the tiling inside the arctic circle.
Henry Cohn, Noam Elkies, and James Propp \cite{CEP} took this step,
and gave a formula governing the behavior of random tilings
inside the circle.
One consequence of their formula
is that regions within the arctic circle
that are macroscopically separated
(i.e., separated by distances whose ratio to $n$ is bounded below)
exhibit different local statistics
under random tiling.
Thus, the extreme homogeneity of a random tiling
in the four regions outside of the arctic circle
is in sharp contrast with the total non-homogeneity
that prevails inside the circle.

Our proof that the boundary between these two domains is circular
is more complicated than we would like.
However, if one is content with the more modest, qualitative goal
of seeing that there must exist some non-homogeneity
in the statistics of the tiling,
then a simple explanation exists.
Recall from \cite{elkies} that domino tilings of the Aztec diamond
are in bijection with certain ``height-functions'' ---
integer-valued functions on the set of lattice points internal to 
the Aztec diamond,
satisfying certain local constraints and boundary conditions.
The local constraints force the height-function
to satisfy a Lipschitz condition,
so that it cannot change too rapidly.
Consider the real-valued function obtained by averaging
all the height-functions that correspond to tilings
(the ``average height-function'').
It too will satisfy a Lipschitz condition.
The difference in average height between two neighboring vertices
is an easily-calculated function of the local tiling statistics,
so that if the local statistics were to be homogeneous,
the average height-function at lattice-points $(i,j)$
would have to be well-approximated by a linear function $ai+bj+c$
for suitable constants $a,b,c$.
However, one can check that the boundary conditions
for Aztec diamond height-functions
are not consistent with any single choice of values for $a,b,c$,
because as one travels around the boundary of the Aztec diamond
the height alternately increases and decreases.

If we were to replace the Aztec diamond by a square
in the preceding analysis,
we would find that the height-functions associated with tilings
are all essentially constant on the boundary of the region,
so that one could take $a=b=0$ and $c$ arbitrary
and get a good approximation along the boundary.
Indeed, Burton and Pemantle's work \cite{burton} shows that
random domino tilings of large squares
are statistically homogeneous away from the boundary
(in a suitable asymptotic sense).
Subsequent work on local statistics for domino tilings 
shows that, in the liquid region, the tiling locally converges 
to an appropriate translation-invariant Gibbs measure 
determined by the macroscopic slope; 
in particular, at the center one obtains 
the symmetric Gibbs measure 
(the Burton–Pemantle/maximal-entropy statistics 
on the infinite square grid);
see \cite{kenyon} and later works.

Shortly after the current article was distributed as a preprint, 
Kenyon \cite{kenyon} proved that 
the statistics at the very center of the Aztec diamond
converge in the limit to Burton-Pemantle statistics.


\begin{ack}
We thank David Aldous and Persi Diaconis for helpful conversations.
Thanks also to Sameera Iyengar, who wrote the first program
for generating random domino tilings by the shuffling algorithm.
We are also grateful to the anonymous referee who made many helpful suggestions.

\
\end{ack}

\begin{funding}
This research was supported by NSF grant DMS 9206374
and NSA grant MDA904-92-H-3060, and by an NSF Graduate Fellowship.
\end{funding}



\begin{thebibliography}{99}

\bibitem{billingsley}
P. Billingsley, \emph{Probability and measure}.
2nd edition, Wiley, 1976

\bibitem{burton}
R.~Burton and R.~Pemantle,
Local characteristics, entropy and limit theorems for spanning
trees and domino tilings via transfer-impedances.
\emph{Ann. Probab.} \textbf{21} (1993), 1329--1371

\bibitem{CEP}
H.~Cohn, N.~Elkies, and J.~Propp,
Local statistics for random domino tilings of the Aztec diamond.
\emph{Duke Math. J.} \textbf{85} (1996), 117--166

\bibitem{CKP}
H.~Cohn, R.~Kenyon, and J.~Propp,
A variational principle for domino tilings.
\emph{J. Amer. Math. Soc.} \textbf{14} (2001), no.~2, 297--346

\bibitem{elkies}
N.~Elkies, G.~Kuperberg, M.~Larsen, and J.~Propp, 
Alternating sign matrices and domino tilings.
\emph{J. Algebr. Comb.} \textbf{1} (1992), 111--132, 219--234 

\bibitem{johansson}
K.~Johansson,
The arctic circle boundary and the Airy process, 
\emph{Ann. Probab.} \textbf{33} (2005), no.~1, 1--30

\bibitem{kasteleyn}
P.W.~Kasteleyn,
The statistics of dimers on a lattice, I. The number of dimer
arrangements on a quadratic lattice.
\emph{Physica} \textbf{27} (1961), 1209--1225

\bibitem{kenyon}
R.~Kenyon,
Local statistics of lattice dimers.
\emph{Ann. Inst. H. Poincaré Probab. Statist.}
\textbf{33} (1997), 591--618

\bibitem{KOS}
R.~Kenyon, A.~Okounkov, and S.~Sheffield,
Dimers and amoebae.
\emph{Ann. of Math.} \textbf{163} (2006), no.~3, 1019--1056

\bibitem{liggett}
T.~Liggett,
\emph{Interacting particle systems}.
Springer-Verlag, 1985

\bibitem{propp}
J.~Propp, Generalized domino-shuffling.
\emph{Theoret. Comput. Sci.} \textbf{303} (2003), 267--301

\bibitem{rost}
H.~Rost, Non-equilibrium behavior of a many-particle system: 
density profile and local equilibria,
\emph{Z. Wahrsch. Verw. Gebiete}
\textbf{58} (1981), 41--53 

\bibitem{smythe}
R.~T.~Smythe and J.~C.~Wierman,
\emph{First-passage percolation on the square lattice}.
Springer-Verlag, 1978

\end{thebibliography}
\end{document}